\documentclass[10pt]{amsart}
\usepackage[T1]{fontenc}
\usepackage{lmodern}
\usepackage{geometry}
\usepackage[latin1] {inputenc}
\usepackage{amsmath}
\usepackage{amsfonts, amssymb, textcomp}
\usepackage[colorlinks=flase, linkcolor=red,urlcolor=green, citecolor=blue]{hyperref}
\usepackage{subeqnarray}

\usepackage{latexsym}
\usepackage{fancyhdr}
\usepackage{longtable}
\usepackage{amsmath, amssymb}
\usepackage{graphicx}

\numberwithin{equation}{section}

\theoremstyle{plain}
\newtheorem{theorem}{Theorem}[section]
\newtheorem{lemma}[theorem]{Lemma}
\newtheorem{proposition}[theorem]{Proposition}
\newtheorem{corollary}[theorem]{Corollary}
\theoremstyle{definition}
\newtheorem{definition}[theorem]{Definition}

\newtheorem{remark}[theorem]{Remark}

\newcommand{\RR}{\mathbb{R}}

\newcommand{\NN}{\mathbb{N}}

\let\on=\operatorname

\newenvironment{demo}[1]{\par\smallskip\noindent{\bf #1.}}{\par\smallskip}

\makeatletter
\@namedef{subjclassname@2020}{%
  \textup{2020} Mathematics Subject Classification}
\makeatother

\title[On equivalent reformulations of condition derivation closedness]
{On equivalent reformulations of condition derivation closedness}

\author[G.~Schindl]{Gerhard Schindl}

\address{G.~Schindl: Fakult\"at f\"ur Mathematik, Universit\"at Wien, Oskar-Morgenstern-Platz~1, A-1090 Wien, Austria.}
\email{gerhard.schindl@univie.ac.at}

\begin{document}

\begin{abstract}
We provide a detailed study of equivalent reformulations of the classical condition derivation closedness for weight sequences. This property is crucial in many applications since it ensures that the corresponding weighted spaces satisfy the desired property of stability under taking derivatives. Moreover, it yields technical estimates which are required in various directions and proofs and so the knowledge that different conditions are equivalent to this crucial requirement is important and useful when dealing with weighted spaces. Finally, we construct an explicit (counter-)example which shows that without the standard condition log-convexity the full characterization fails in general and comment on the mixed setting as well.
\end{abstract}

\thanks{This research was funded in whole or in part by the Austrian Science Fund (FWF) 10.55776/PAT9445424}
\date{\today}

\maketitle

\section{Introduction}
When working with weighted function spaces for which the defining weight is expressed in terms of a \emph{weight sequence} $\mathbf{M}=(M_j)_{j\in\NN}$ (we write $\NN:=\{0,1,2,\dots\}$ and $\NN_{>0}:=\NN\setminus\{0\}$) resp. in terms of the corresponding \emph{associated weight function} $\omega_{\mathbf{M}}$ (see \eqref{assofunc}), then growth and regularity properties for $\mathbf{M}$ are unavoidable and natural. A classical and very prominent one is the condition $(M.2)'$ resp. \emph{stability under differential operators} in the notation of H. Komatsu (see \cite{Komatsu73}) and which reads as follows:
\begin{equation}\label{dckomatsu}
\exists\;A,H>0\;\forall\;j\in\NN:\;\;\;M_{j+1}\le AH^jM_j.
\end{equation}
Obviously, this is equivalent to the more compact form
\begin{equation}\label{dc}
\exists\;D\ge 1\;\forall\;j\in\mathbb{N}:\;\;\;M_{j+1}\le D^{j+1}M_j,
\end{equation}
resp. write $D^j$ instead and focus on all $j\in\NN_{>0}$. We denote this condition by \emph{derivation closedness} and abbreviate it by \hypertarget{dc}{$(\on{dc})$}. The first aim of this article is to provide an exhaustive study of equivalent reformulations of this condition in terms of growth conditions for $\mathbf{M}$, the associated function $\omega_{\mathbf{M}}$ (see \eqref{assofunc}) and the corresponding \emph{counting function} $\Sigma_{\mathbf{M}}$ (see \eqref{counting}). For this we need that $\mathbf{M}$ satisfies only elementary growth conditions; most importantly it turns out that one requires as a basic property that $\mathbf{M}$ is \emph{log-convex;} i.e.
$$\forall\;j\in\NN_{>0}:\;\;\;M_j^2\le M_{j+1}M_{j-1},$$
and this is precisely (M.1) in \cite{Komatsu73}. For applications in proofs this knowledge is crucial and convenient since occasionally the different variants of \hyperlink{dc}{$(\on{dc})$} are required naturally. We also comment on the mixed setting and hence this work is in the spirit of the author's recent papers \cite{modgrowthstrange} and \cite{modgrowthstrangeII} in which the analogous goal has been treated concerning the classical and stronger condition \emph{moderate growth (stability under ultradifferential operators)}
\begin{equation*}\label{mg}
\exists\;C\ge 1\;\forall\;j,k\in\NN:\;\;\;M_{j+k}\le C^{j+k+1}M_jM_k.
\end{equation*}
This corresponds to $(M.2)$ in \cite{Komatsu73} and we write \hypertarget{mg}{$(\on{mg})$} for it. More precisely, the motivation behind this paper is the following:

\begin{itemize}
\item[$(A)$] Several equivalent reformulations of \hyperlink{dc}{$(\on{dc})$} can already be found in the founding book by S. Mandelbrojt \cite{mandelbrojtbook}; more precisely confirm 6.9 IV on p. 248 \& 249, (6.9.4) resp. (6.9.5) on p. 246; 1.8.VIII on p. 24 and finally the main result Thm. 1.6.II on p. 14. The proof of this statement is contained in \cite[Ch. I, Sect. 1.6, p. 9-14]{mandelbrojtbook} and in fact it deals with a much more general notation involving the regularized sequence w.r.t. a \emph{regularizing function} $\omega$. In the standard setting dealing with the \emph{log-convex minorant} $\mathbf{M}^{\on{lc}}$ this corresponds to $\omega\equiv+\infty$ (``r\'{e}gularisation convexe'') and we also refer to \cite{regularnew} for more details and explanations.

    Thus, if one considers \emph{log-convex sequences} then we have $\mathbf{M}=\mathbf{M}^{\on{lc}}$ and $\omega\equiv+\infty$, the requirement $\liminf_{t\rightarrow+\infty}\frac{\omega(t)}{t}>0$ in \cite[Thm. 1.6.II]{mandelbrojtbook} (see also item $d.$ on p. 10) is trivial and the \emph{trace function} $A$ is continuous (cf. \cite[Sect. 1.3, p. 6]{mandelbrojtbook} and \cite{regularnew}). Therefore, by Mandelbrojt's approach for log-convex sequences a full characterization is available and note that in this setting it is also natural to assume that
    \begin{equation}\label{rootinfinity}
    \lim_{j\rightarrow+\infty}(M_j)^{\frac{1}{j}}=+\infty;
    \end{equation}
    see again \cite[Ch. I]{mandelbrojtbook} and \cite{regularnew}.

\item[$(B)$] We have also found that in a work by W. Matsumoto, see \cite[Thm. 5]{matsumoto}, an exhaustive characterization for \hyperlink{dc}{$(\on{dc})$} is mentioned and rephrasing the conditions from \cite{mandelbrojtbook} but without giving a precise proof and referring to S. Mandelbrojt's book on \cite[p. 669]{matsumoto}. In \cite{matsumoto} basically it was assumed that $\mathbf{M}$ is \emph{log-convex} and satisfies \eqref{rootinfinity}; see p. 668. Note that condition $(D)$ on \cite[p. 669]{matsumoto} should correspond to \eqref{dc} but there is a misprint since instead of $C$ one should read $C^{j+1}$: the stated condition implies by iteration $\sup_{j\in\NN_{>0}}(M_j)^{\frac{1}{j}}<+\infty$ and hence is contradicting \eqref{rootinfinity}.

    The main content of \cite{matsumoto}, dealing with different reformulations of \hyperlink{mg}{$(\on{mg})$} (in \cite{matsumoto} this is called \emph{separativity condition}), admits an interesting observation: In the main result \cite[Thm. 1]{matsumoto} the exhaustive characterization of \hyperlink{mg}{$(\on{mg})$} is established, and in \cite[Thm. 2]{matsumoto} the equivalence of several conditions which are necessary for having \hyperlink{mg}{$(\on{mg})$} is studied. Condition \cite[Thm. 2, i)]{matsumoto} means that each $\mathbf{M}$ satisfying \hyperlink{mg}{$(\on{mg})$} is growing at most like $j\mapsto j!^s$ for some $s>0$; i.e. $\mathbf{M}$ can be dominated by some \emph{Gevrey sequence} $\mathbf{G}^s:=(j!^s)_{j\in\NN}$ of index/order $s>0$. However, the converse is in general not true; a (counter-)example is constructed in \cite[p. 672-673]{matsumoto}.\vspace{6pt}

    On the other hand from \cite[Thm. 1.6.II]{mandelbrojtbook} (and \cite[Thm. 5]{matsumoto}) it does follow that \hyperlink{dc}{$(\on{dc})$} is equivalent to
    \begin{equation}\label{dcMandelbrojt}
\sup_{j\in\NN_{>0}}\frac{\log(M_j)}{j^2}<+\infty\Leftrightarrow\sup_{j\in\NN_{>0}}\frac{\log((M_j)^{\frac{1}{j}})}{j}<+\infty.
\end{equation}
    In view of all these results and observations one concludes that a qualitative difference between \hyperlink{dc}{$(\on{dc})$} and \hyperlink{mg}{$(\on{mg})$} appears: Whereas \hyperlink{mg}{$(\on{mg})$} does imply only necessarily an upper bound for the growth of $\mathbf{M}$ expressed in terms of some Gevrey sequence, for \hyperlink{dc}{$(\on{dc})$} condition \eqref{dcMandelbrojt} provides an upper bound expressed in terms of some \emph{$q$-Gevrey sequence} $(\ell^{j^2})_{j\in\NN}$, $\ell>1$, and this growth restriction does already characterize \hyperlink{dc}{$(\on{dc})$}.\vspace{6pt}

    First, we believed that there might be a gap in the reformulation of \hyperlink{dc}{$(\on{dc})$} (and hence in the proof of \cite{mandelbrojtbook}) since from this point of view the difference between \hyperlink{dc}{$(\on{dc})$} and \hyperlink{mg}{$(\on{mg})$} is not clear and one might think to transfer the construction of the (counter-)example from Matsumoto to condition \hyperlink{dc}{$(\on{dc})$}. However, we have recovered the proof, see the main result Theorem \ref{maindcthm}, and it turned out that the characterization for \hyperlink{dc}{$(\on{dc})$} is correct.\vspace{6pt}

    Based on this study, and this is another interesting phenomenon in this context, in the second main statement Proposition \ref{Mcounterexample} we have been able to construct a technical \emph{(counter-)example} which shows that for the full characterization of \hyperlink{dc}{$(\on{dc})$} (in Theorem \ref{maindcthm}) the basic condition \emph{log-convexity is unavoidable.} This is consistent with the assumptions on $\mathbf{M}$ resp. the regularizing function $\omega$ in \cite[Thm. 1.6.II]{mandelbrojtbook} (and in \cite[Thm. 5]{matsumoto}) and illustrates the importance of the precise growth of the regularizing function in \cite{mandelbrojtbook}. In view of this observation, see in particular $(e)$ in Proposition \ref{Mcounterexample}.

   However, the constructed sequence does still have ``nice'' regularity properties: $j\mapsto(M_j)^{\frac{1}{j}}$ is non-decreasing and, moreover, \eqref{rootinfinity} is also valid. Occasionally these are standard assumptions in the ultradifferentiable weight sequence setting; see the approach by H\"{o}rmander e.g. in \cite{hoermander}. Here, the spaces are defined in terms of the sequence of roots $L_j:=(M_{j})^{\frac{1}{j}}$, $j\in\NN_{>0}$, and set $L_0:=1$. Crucial conditions are expressed in terms of this modification and hence instead of \hyperlink{dc}{$(\on{dc})$} it is required that
     \begin{equation}\label{rootaincr}
\exists\;A\ge 1\;\forall\;j\in\NN_{>0}:\;\;\;(M_{j+1})^{\frac{1}{j+1}}\le A(M_{j})^{\frac{1}{j}}.
\end{equation}
Roughly speaking our constructed sequence illustrates the fact that within the ``H\"{o}rmander setting'' in general we cannot expect the full equivalence of all conditions: there exist $\mathbf{M}$ such that \eqref{dcMandelbrojt} holds but both \hyperlink{dc}{$(\on{dc})$} and \eqref{rootaincr} fail!

\item[$(C)$] Another motivation for this article was coming from the author's recent works \cite{ultradifferentiablecomparison} and \cite{convolvedmatrix}.

    For concrete applications and when dealing with technical issues in proofs it is relevant and unavoidable to have knowledge on the precise growth of $t\mapsto\omega_{\mathbf{M}}(h't)-\omega_{\mathbf{M}}(ht)$ as $t\rightarrow+\infty$ with (different) real parameters $h'>h>0$. The motivation behind this question is to get ``enough (sufficient) growth'' when $h'>h$ accordingly and one can think here in two directions: if $h>0$ is given, then choose $h'>h$ resp. when $h'>0$ is given then select $h<h'$ in a proper way in order to have the required growth.\vspace{6pt}

    This question becomes important when dealing with \emph{dilatation-type} weight function systems/matrices in the terminology of \cite{weightedentireinclusion1}, \cite{weightedentireinclusion2} and \cite{ultradifferentiablecomparison} since here the parameters $h,h'>0$ play the crucial role of the \emph{matrix index.} In \cite[Lemma 2.9]{weightedentireinclusion1} we have shown that for any log-convex sequence $\mathbf{M}=(M_j)_{j\in\NN}$ satisfying \eqref{rootinfinity} (and $M_0=1$, cf. \cite[Def. 2.4]{weightedentireinclusion1}) we get
\begin{equation}\label{L29weightentire}
\forall\;h'>h>0:\;\;\;\lim_{t\rightarrow+\infty}\omega_{\mathbf{M}}(h't)-\omega_{\mathbf{M}}(ht)=+\infty.
\end{equation}
Note that normalization $M_0=1$ is not effecting the proof of this result (it crucially involves \eqref{assointrepr}) and \eqref{L29weightentire} is sufficient and appropriate in order to treat weighted $L^{\infty}$-type norms as it is done in the weighted entire setting in \cite{weightedentireinclusion1} and \cite{weightedentireinclusion2} and for weighted $L^{\infty}$-type spaces of functions $f: \mathbb{R}^d\rightarrow\mathbb{C}$ in \cite{ultradifferentiablecomparison}. However, when dealing with weighted $L^p$-type norms with $1\le p<+\infty$, cf. the approach in \cite{ultradifferentiablecomparison}, then one is interested in having boundedness of the integrals $\int_{\RR^d}e^{\omega_{\mathbf{M}}(ht)-\omega_{\mathbf{M}}(h't)}dt$ resp. $\int_{\RR^d}e^{p\omega_{\mathbf{M}}(ht)-p\omega_{\mathbf{M}}(h't)}dt$ (usually, in the literature $\omega_{\mathbf{M}}$ is in this case radially extended to whole $\mathbb{R}^d$). And for $p=1$ this question becomes also relevant in ultradifferentiable Paley-Wiener-type results; see e.g. the proofs of \cite[Thm. 9.1]{Komatsu73} and of \cite[Thm. 4.5]{convolvedmatrix}.\vspace{6pt}

It turns out that \hyperlink{dc}{$(\on{dc})$} is unavoidable for having finiteness/boundedness of the aforementioned integrals; indeed \hyperlink{dc}{$(\on{dc})$} guarantees finiteness of the corresponding integral(s) \emph{uniformly for any $p\in[1,+\infty)$.}
\end{itemize}

This paper is structured as follows: In Section \ref{weightsequencesection} we summarize more definitions and growth properties concerning weight sequences and associated weight functions. In Section \ref{mainresultsection} we prove the main characterizing result Theorem \ref{maindcthm} and split it into several partial results in order to give a clear picture and to illustrate the connections between the equivalent conditions; see Theorems \ref{firstcharprop}, \ref{secondcharprop}, \ref{fourthcharprop}, \ref{spacedc}, and Proposition \ref{thirdcharprop}. We are able to present an exhaustive list of equivalent reformulations of \hyperlink{dc}{$(\on{dc})$} and extend the known conditions by \eqref{newseriescond} and \eqref{newintcond}; the latter one is tailor-made for the questions discussed in $(C)$. Note that for some parts of the proof log-convexity resp. \eqref{rootinfinity} is \emph{not required necessarily} (cf. Theorem \ref{firstcharprop}) and for some parts we apply recent weight matrix techniques; see the proof of Theorem \ref{fourthcharprop}. Section \ref{counterexsection} is dedicated to the second main statement in this work: the construction of the (counter-)example in Proposition \ref{Mcounterexample}. Finally, in Section \ref{supplementarysection} we discuss some additional issues: In Section \ref{stabilitysection} we comment on the stability of \hyperlink{dc}{$(\on{dc})$} under several natural weight sequence constructions. And in Section \ref{mixedsection} we study which parts of Theorem \ref{maindcthm} transfer to the mixed setting. The mixed notion is becoming relevant and natural for the general weight matrix setting and, in particular, for weight functions in the sense of Braun-Meise-Taylor by involving the \emph{associated weight matrix} $\mathcal{M}_{\omega}$; see \cite[Sect. 4.1 \& 5]{compositionpaper}. However, as it turns out only some parts can be recovered and this is the reason why not formulating the central statement even in terms of mixed conditions.

\section{Weight sequences}\label{weightsequencesection}
We introduce more relevant notation for sequences. $\mathbf{M}$ is called \emph{normalized} if $1=M_0\le M_1$. The corresponding sequence of quotients $\mu=(\mu_j)_{j\in\NN_{>0}}$ is given by $\mu_j:=\frac{M_j}{M_{j-1}}$ for $j\in\NN_{>0}$ and set $\mu_0:=1$. Analogous Greek letters will used then considering other sequences; e.g. write $\nu$ for the quotients of the sequence $\mathbf{N}$. Log-convexity precisely means that $j\mapsto\mu_j$, $j\in\NN_{>0}$, is non-decreasing. On the other hand, recall that for any log-convex sequence $\mathbf{M}$ we have that $j\mapsto\left(\frac{M_j}{M_0}\right)^{\frac{1}{j}}$ is non-decreasing: apply (the proof of) \cite[Lemma 2.0.4]{diploma} to the sequence $\widetilde{\mathbf{M}}:=(M_j/M_0)_{j\in\NN}$. Therefore, when $M_0=1$ and which can be assumed w.l.o.g. when being interested in the corresponding weighted spaces by considering $\widetilde{\mathbf{M}}$ instead, we have that log-convexity is stronger than $j\mapsto(M_j)^{\frac{1}{j}}$ being non-decreasing. This difference is also emphasized by the (counter-)example in Proposition \ref{Mcounterexample}.

For given $\mathbf{M}\in\RR_{>0}^{\NN}$ define $$\mathbf{M}_{\iota}:=\liminf_{j\rightarrow+\infty}\left(\frac{M_j}{M_0}\right)^{\frac{1}{j}},\hspace{15pt}\mathbf{M}_{\sigma}:=\limsup_{j\rightarrow+\infty}\left(\frac{M_j}{M_0}\right)^{\frac{1}{j}};$$
see \cite[p. 180, $(2.3)$ \& p. 193]{regularnew}. Note that $\mathbf{M}_{\iota}>0$ if and only if $\inf_{j\in\NN_{>0}}(M_j)^{\frac{1}{j}}>0$. By \cite[Sect. 2.2, Lemma 2.1]{regularnew} (and \cite[Rem. 4.1]{conjugateweightfunction}) it holds that for any log-convex $\mathbf{M}$ one has $\mathbf{M}_{\iota}=\mathbf{M}_{\sigma}=\lim_{j\rightarrow+\infty}(M_j)^{\frac{1}{j}}=\lim_{j\rightarrow+\infty}\mu_j\in(0,+\infty]$. Occasionally we consider the following (sub)set of weight sequences:
$$\hypertarget{LCset}{\mathcal{LC}}:=\{\mathbf{M}\in\RR_{>0}^{\NN}:\;\mathbf{M}\;\text{is normalized, log-convex},\;\lim_{j\rightarrow+\infty}(M_j)^{\frac{1}{j}}=+\infty\}.$$
Thus $\mathbf{M}\in\hyperlink{LCset}{\mathcal{LC}}$ if and only if $1=\mu_0\le\mu_1\le\dots$ and $\lim_{j\rightarrow+\infty}\mu_j=+\infty$ (see e.g. \cite[p. 104]{compositionpaper}) and there is a one-to-one correspondence between $\mathbf{M}$ and $\mu=(\mu_j)_j$ by taking $M_j:=\prod_{k=0}^j\mu_k$.\vspace{6pt}

For the following definition we refer to \cite[Chapitre I]{mandelbrojtbook}, \cite[Sect. 3]{Komatsu73} and the more recent work \cite{regularnew} in which additional properties and explanations are summarized and non-standard cases are treated. Let $\mathbf{M}=(M_j)_{j\in\NN}\in\RR_{>0}^{\NN}$ be given, then the \emph{associated (weight) function} $\omega_{\mathbf{M}}: [0,+\infty)\rightarrow[0,+\infty)\cup\{+\infty\}$ is defined as follows:
\begin{equation}\label{assofunc}
\omega_{\mathbf{M}}(t):=\sup_{j\in\NN}\log\left(\frac{M_0t^j}{M_j}\right),\qquad t\ge 0,
\end{equation}
with the conventions $0^0:=1$ and $\log(0)=-\infty$. This ensures $\omega_{\mathbf{M}}(0)=0$ and $\omega_{\mathbf{M}}(t)\ge 0$ for any $t\ge 0$ since $\frac{t^0M_0}{M_0}=1$ for all $t\ge 0$. \eqref{assofunc} corresponds to \cite[$(3.1)$]{Komatsu73} and we immediately have that $\omega_{\mathbf{M}}$ is non-decreasing and satisfying $\lim_{t\rightarrow+\infty}\omega_{\mathbf{M}}(t)=+\infty$. However, the fact $\omega_{\mathbf{M}}: [0,+\infty)\rightarrow[0,+\infty)$ (well definedness) is equivalent to assume $\lim_{j\rightarrow+\infty}(M_j)^{\frac{1}{j}}=+\infty$ (see \cite[Lemma 2.2]{regularnew}).

Let $\mathbf{M}$ be a log-convex sequence satisfying $\lim_{j\rightarrow+\infty}(M_j)^{\frac{1}{j}}=+\infty$, hence $\lim_{j\rightarrow+\infty}\mu_j=+\infty$ and define the \emph{counting function}
\begin{equation}\label{counting}
\Sigma_{\mathbf{M}}(t):=|\{j\in\NN_{>0}:\;\;\;\mu_j\le t\}|.
\end{equation}
By the known \emph{integral representation formula,} see \cite[1.8. III]{mandelbrojtbook}, \cite[$(3.11)$]{Komatsu73} and also \cite[Lemma 2.5]{regularnew}, for all such sequences one has
\begin{equation}\label{assointrepr}
\forall\;t\ge 0:\;\;\;\omega_{\mathbf{M}}(t)=\int_0^t\frac{\Sigma_{\mathbf{M}}(u)}{u}du=\int_{\mu_1}^t\frac{\Sigma_{\mathbf{M}}(u)}{u}du.
\end{equation}
Thus the integral near $0$ is vanishing since $\Sigma_{\mathbf{M}}(t)=0$ for $t\in[0,\mu_1)$ (empty set has cardinality $0$). Finally, let $\mathbf{M}$ be log-convex and such that $\lim_{j\rightarrow+\infty}(M_j)^{\frac{1}{j}}=+\infty$, then one can compute $\mathbf{M}$ by involving $\omega_{\mathbf{M}}$ as follows, see \cite[Chapitre I, 1.4, 1.8]{mandelbrojtbook}, \cite[Prop. 3.2]{Komatsu73} and more recently \cite{regularnew}:
\begin{equation}\label{Prop32Komatsu}
M_j=M_0\sup_{t\ge 0}\frac{t^j}{\exp(\omega_{M}(t))},\;\;\;j\in\NN.
\end{equation}
If $\mathbf{M}$ satisfies $\lim_{j\rightarrow+\infty}(M_j)^{\frac{1}{j}}=+\infty$ but is not log-convex, then the right-hand side of \eqref{Prop32Komatsu} yields $M^{\on{lc}}_j$ with $\mathbf{M}^{\on{lc}}$ denoting the \emph{log-convex minorant} of $\mathbf{M}$.\vspace{6pt}

When $\mathbf{M},\mathbf{N}\in\RR_{>0}^{\NN}$ are given, then let us write $\mathbf{M}\hypertarget{preceq}{\preceq}\mathbf{N}$ if
$$\sup_{j\in\NN_{>0}}\left(\frac{M_j}{N_j}\right)^{\frac{1}{j}}<+\infty;$$
i.e. if $M_j\le C^{j+1}N_j$ for some $C\ge 1$ and all $j\in\NN$. This relation is becoming relevant for the inclusion relations in the corresponding weighted settings, for the ultradifferentiable classes see e.g. \cite[Sect. 2.7, Prop. 2.12]{compositionpaper}. We say that $\mathbf{M}$, $\mathbf{N}$ are equivalent, written $\mathbf{M}\hypertarget{approx}{\approx}\mathbf{N}$, if $\mathbf{M}\hyperlink{preceq}{\preceq}\mathbf{N}$ and $\mathbf{N}\hyperlink{preceq}{\preceq}\mathbf{M}$. Obviously, $\widetilde{\mathbf{M}}\hyperlink{approx}{\approx}\mathbf{M}$ with $\widetilde{\mathbf{M}}:=(M_j/M_0)_{j\in\NN}$ and \hyperlink{dc}{$(\on{dc})$}, \eqref{dcMandelbrojt} and \eqref{rootaincr} are clearly preserved under equivalence.

This last comment implies an immediate technical observation:

\begin{lemma}\label{dctrivial}
Let $\mathbf{M}=(M_j)_{j\in\NN}\in\RR_{>0}^{\NN}$ be given, then \hyperlink{dc}{$(\on{dc})$}, \eqref{dcMandelbrojt} and \eqref{rootaincr} are valid if
\begin{itemize}
\item[$(i)$] $\mathbf{M}$ satisfies
\begin{equation}\label{dctrivial}
0<\mathbf{M}_{\iota}\le\mathbf{M}_{\sigma}<+\infty.
\end{equation}
\item[$(ii)$] A particular case for which \eqref{dctrivial} holds is when $\mathbf{M}$ is log-convex and satisfies $\lim_{j\rightarrow+\infty}(M_j)^{\frac{1}{j}}<+\infty$.
\end{itemize}
\end{lemma}

\emph{Note:} $(ii)$ yields the fact that in the main results for which log-convexity of $\mathbf{M}$ is required (cf. Theorem \ref{maindcthm}) the (additional) assumption $(\mathbf{M}_{\iota}=)\lim_{j\rightarrow+\infty}(M_j)^{\frac{1}{j}}=+\infty$ is natural and provides non-trivial information.

\demo{Proof}
$(i)$ \eqref{dctrivial} does precisely mean
$$\exists\;0<c<1<C<+\infty\;\forall\;j\in\NN:\;\;\;c^{j+1}\le M_j\le C^{j+1};$$
i.e. $\mathbf{M}$ and the constant sequence $\mathbf{1}:=(1_j)_{j\in\NN}$ are equivalent. Since $\mathbf{1}$ obviously enjoys the desired properties and since they are preserved under equivalence we are done.

This should be compared with the remark just after \cite[Cor. 4.1, p. 193]{regularnew} and, indeed, for the arguments there one can also allow \eqref{dctrivial} and having $\mathbf{M}_{\iota}=\mathbf{M}_{\sigma}$ as well.\vspace{6pt}

$(ii)$ In this case, as explained before, one has $0<\mathbf{M}_{\iota}=\mathbf{M}_{\sigma}=\lim_{j\rightarrow+\infty}(M_j)^{\frac{1}{j}}<+\infty$ and hence \eqref{dctrivial} follows.
\qed\enddemo

\section{Main results}\label{mainresultsection}
We divide the characterization into several parts in order to make the characterization more clear since some parts need different assumptions on $\mathbf{M}$ resp. different techniques are required.

\subsection{First statement}\label{firstmainresult}
Our first main result is the following characterization:

\begin{theorem}\label{firstcharprop}
Let $\mathbf{M}=(M_j)_{j\in\NN}\in\RR_{>0}^{\NN}$ be given and let $p\in[1,+\infty)$ be fixed. Consider the following assertions:
\begin{itemize}
\item[$(i)$] $\mathbf{M}$ satisfies \hyperlink{dc}{$(\on{dc})$} (see \eqref{dc}); i.e.
$$\exists\;D\ge 1\;\forall\;j\in\mathbb{N}:\;\;\;M_{j+1}\le D^{j+1}M_j.$$

\item[$(ii)$] The sequence of quotients $\mu=(\mu_j)_{j\in\NN_{>0}}$ satisfies
\begin{equation}\label{mudc}
\exists\;D\ge 1\;\forall\;j\in\mathbb{N}:\;\;\;\mu_{j+1}\le D^{j+1}.
\end{equation}

\item[$(iii)$] The sequence of quotients $\mu=(\mu_j)_{j\in\NN_{>0}}$ satisfies
\begin{equation}\label{mudiff}
\exists\;A\ge 1\;\forall\;j\in\mathbb{N}_{>0}:\;\;\;\mu_{j+1}-\mu_j\le A^j.
\end{equation}

\item[$(iv)$] We get
\begin{equation}\label{newseriescond}
\exists\;h'>h>0:\;\;\;\sum_{j=1}^{+\infty}\left(\frac{h}{h'}\right)^{pj}(\mu_{j+1}-\mu_j)<+\infty.
\end{equation}

\item[$(v)$] $\mathbf{M}$ satisfies \eqref{rootaincr}; i.e.
$$\exists\;A\ge 1\;\forall\;j\in\NN_{>0}:\;\;\;(M_{j+1})^{\frac{1}{j+1}}\le A(M_{j})^{\frac{1}{j}},$$
which is equivalent to
$$\forall\;d\in\NN_{>0}\;\exists\;A\ge 1\;\forall\;j\in\NN_{>0}:\;\;\;(M_{j+d})^{\frac{1}{j+d}}\le A(M_{j})^{\frac{1}{j}}.$$

\item[$(vi)$] $\mathbf{M}$ satisfies \eqref{dcMandelbrojt}; i.e.
$$\exists\;C\ge 1\;\forall\;j\in\NN:\;\;\;M_j\le C^{j^2+1}.$$
\end{itemize}

Then $(v)\Rightarrow(i)\Leftrightarrow(ii)\Leftrightarrow(iii)\Rightarrow(iv)$ and $(v)\Rightarrow(vi)$ is valid. When in addition $\mathbf{M}_{\iota}>0$, then also $(i)\Rightarrow(v)$ and when $\mathbf{M}$ is log-convex, then $(i)\Rightarrow(v)$ and $(iv)\Rightarrow(iii)$. Thus for log-convex $\mathbf{M}$ the first five assertions are equivalent.
\end{theorem}

\emph{Note:}

\begin{itemize}
\item[$(*)$] The stated equivalence in assertion $(v)$ is trivial by iteration.

\item[$(*)$] If $\mathbf{M}$ is log-convex and such that $\mathbf{M}_{\iota}\in(0,+\infty)$, then this result is consistent with Lemma \ref{dctrivial}. However, in any case for this theorem it is \emph{not required necessarily that} $\lim_{j\rightarrow+\infty}(M_j)^{\frac{1}{j}}=+\infty$ but, on the other hand, this condition clearly implies $\mathbf{M}_{\iota}>0$.

\item[$(*)$] As mentioned in the introduction, occasionally in the literature in \eqref{dc} one can read $AB^j$ instead of $D^{j+1}$ but which is obviously equivalent (set $D:=\max\{A,B\}$). Similarly also in $(vi)$ one could replace $C^{j^2+1}$ by $A^{j^2}B$ for some constants $A,B\ge 1$.

\item[$(*)$] The precise value of $p\in[1,+\infty)$ in \eqref{newseriescond} is not effecting the validity of the above characterization and hence one has uniformity w.r.t. this parameter. As it turns out in the proof, only the precise relation between $h'$ and $h$ is depending on $p$.

\item[$(*)$] In \eqref{newseriescond} it is natural to restrict w.l.o.g. to $h'>h$. Consider $0<h'\le h$ and note: When $\lim_{j\rightarrow+\infty}\mu_j=+\infty$, then the sum is divergent since the first factor in each summand is at least $1$ and $\sum_{j=1}^{+\infty}(\mu_{j+1}-\mu_j)=+\infty$ because this corresponds to a telescopic series. If $\inf_{j\in\NN_{>0}}\mu_{j}-\mu_{j}>0$ then the series in \eqref{newseriescond} fails to be convergent too and when $j\mapsto\mu_j$ is non-decreasing, i.e. $\mathbf{M}$ is log-convex, then it suffices to have $\inf_{k\in\NN_{>0}}\mu_{j_k+1}-\mu_{j_k}>0$ for a subsequence $(j_k)_{k\in\NN_{>0}}$. Again, these properties are not depending on the value of $p$.
\end{itemize}

\demo{Proof}
Before showing the implications recall that log-convexity implies $\mathbf{M}_{\iota}>0$.\vspace{6pt}

$(v)\Rightarrow(i)$ First, $(M_{j+1})^{\frac{1}{j+1}}\le A(M_{j})^{\frac{1}{j}}\Leftrightarrow M_{j+1}\le A^{j+1}M_j(M_j)^{\frac{1}{j}}$ for all $j\in\NN_{>0}$. Then, by iteration of \eqref{rootaincr} we get $(M_{j+1})^{\frac{1}{j+1}}\le A(M_{j})^{\frac{1}{j}}\le A^2(M_{j-1})^{\frac{1}{j-1}}\le\dots\le A^jM_1$ and hence $M_{j+1}\le A^{j(j+1)}(M_1)^{j+1}$ for all $j\in\NN$. Summarizing,
$$M_{j+1}\le A^{j+1}M_j(M_j)^{\frac{1}{j}}\le A^{j+1}M_jA^{j-1}M_1\le M_1A^{2j}M_j\le(\max\{A^2,M_1\})^{j+1}M_j$$
for all $j\in\NN_{>0}$ and which proves \eqref{dc}. Moreover, the estimate obtained by the iteration shows \eqref{dcMandelbrojt} and hence $(v)\Rightarrow(vi)$, too.\vspace{6pt}

$(i)\Leftrightarrow(ii)$ This is clear by definition of the sequence $\mu$ with the same choice for $D$.\vspace{6pt}

$(ii)\Rightarrow(iii)$ Since $\mu_j\ge 0$ for all $j$ we infer the estimate in $(iii)$ with $A:=D^2$.\vspace{6pt}

$(iii)\Rightarrow(ii)$ W.l.o.g. assume $A\ge 2$ and by iteration we get
$$\forall\;j\in\NN_{>0}:\;\;\;\mu_{j+1}\le A^j+\mu_j\le A^j+A^{j-1}+\mu_{j-1}\le\dots\le\sum_{k=1}^{j}A^k+\mu_1,$$
therefore,
$$\mu_{j+1}\le\sum_{k=1}^{j}A^k+\mu_1=\frac{A^{j+1}-A}{A-1}+\mu_1\le A^{j+1}+\mu_1\le 2(\max\{A,\mu_1\})^{j+1}\le(2\max\{A,\mu_1\})^{j+1}.$$ And since this estimate also holds for $j=0$ we infer $(ii)$ with $D:=2\max\{A,\mu_1\}$.\vspace{6pt}

$(iii)\Rightarrow(iv)$ This is clear for all $h'>A^{1/p}h$ with $A\ge 1$ denoting the constant from $(iii)$.\vspace{6pt}

$(iv)\Rightarrow(iii)$ Assume that assertion $(iii)$ fails, then
$$\forall\;n\in\mathbb{N}_{>0}\;\exists\;j_n\in\mathbb{N}_{>0}:\;\;\;\mu_{j_n+1}-\mu_{j_n}\ge n^{j_n},$$
and hence for all $h'>h>0$
$$\sum_{j=1}^{+\infty}\left(\frac{h}{h'}\right)^{pj}(\mu_{j+1}-\mu_j)\ge\sum_{n\in\mathbb{N}_{>0}}\left(\frac{h}{h'}\right)^{pj_n}(\mu_{j_n+1}-\mu_{j_n})\ge\sum_{n\in\mathbb{N}_{>0}}\left(\frac{hn^{1/p}}{h'}\right)^{pj_n}=+\infty,$$
since for any possible choice $h'>h>0$ there exists $n_0\in\mathbb{N}$ such that $hn^{1/p}>h'$ for all $n\ge n_0$. (Similarly, for any $0<h'\le h$ the series is clearly divergent, too.) However, for the first estimate it is crucial to have log-convexity in order to guarantee that each summand in the series is non-negative ($\mu_{j+1}-\mu_j\ge 0$ holds for all $j$).\vspace{6pt}

$(i)\Rightarrow(v)$ When $\mathbf{M}_{\iota}>0$, i.e. $\inf_{j\in\NN_{>0}}(M_j)^{\frac{1}{j}}>0$ and so $(M_j)^{\frac{1}{j}}\ge\epsilon$ for some $\epsilon\in(0,1]$ and all $j\in\NN_{>0}$, then
$$M_{j+1}\le D^{j+1}M_j\le D^{j+1}M_j\frac{(M_j)^{\frac{1}{j}}}{\epsilon}\le(D\epsilon^{-1})^{j+1}M_j(M_j)^{\frac{1}{j}}$$
for some $D\ge 1$ and all $j\in\NN_{>0}$. This is precisely \eqref{rootaincr} with $A:=D\epsilon^{-1}$ (see step $(v)\Rightarrow(i)$).
\qed\enddemo

\begin{remark}\label{mgcharrem}
Let $\mathbf{M}\in\hyperlink{LCset}{\mathcal{LC}}$. An exhaustive list of equivalent reformulations of condition \hyperlink{mg}{$(\on{mg})$} for $\mathbf{M}$ is given in \cite[Thm. 3.1]{modgrowthstrange}, see also the citations there, and for the study of the mixed conditions in the weight matrix setting the following requirement was becoming relevant in \cite{modgrowthstrange}, \cite{modgrowthstrangeII}:
\begin{equation}\label{mgremequ}
\exists\;d\in\NN_{>0}\;\exists\;B\ge 1\;\forall\;j\in\NN_{>0}:\;\;\;\mu_j\le B(M_{dj})^{\frac{1}{dj}}.
\end{equation}
This condition is generalizing one of the equivalent reformulations: \hyperlink{mg}{$(\on{mg})$} is equivalent to \eqref{mgremequ} with $d=1$ and which is stronger under the general assumption on $\mathbf{M}$ since $j\mapsto(M_j)^{\frac{1}{j}}$ is non-decreasing. It is known that \hyperlink{mg}{$(\on{mg})$} is also equivalent to having the estimate in \eqref{mgremequ} with $d=1$ and $\mu_{j+1}$ instead of $\mu_j$ (see e.g. \cite[Lemma 2.2, p. 4]{whitneyextensionweightmatrix}) and one ask if also in \eqref{mgremequ} one can replace $\mu_j$ by $\mu_{j+1}$.

Indeed, this holds for any log-convex $\mathbf{M}$ satisfying \hyperlink{dc}{$(\on{dc})$}: Log-convexity implies $\mu_{j+1}\ge\mu_j$ for all $j$ and so \eqref{mgremequ} follows. Conversely, by $(i)\Leftrightarrow(v)$ in Theorem \ref{firstcharprop} we can iterate \eqref{rootaincr} $d$-times and get $\mu_{j+1}\le B(M_{dj+d})^{\frac{1}{dj+d}}\le BA^d(M_{dj})^{\frac{1}{dj}}$ for all $j\in\NN_{>0}$; thus $\mu_j$ replaced by $\mu_{j+1}$ in \eqref{mgremequ} holds with $B':=BA^d$ and note that the crucial parameter $d$ is preserved in this equivalence. Recall that there exist sequences satisfying \eqref{mgremequ} for $d\ge 2$, violating \hyperlink{mg}{$(\on{mg})$} but having \hyperlink{dc}{$(\on{dc})$}; see \cite[Rem. 4.5 $(ii)$, p. 23]{modgrowthstrange}: Indeed, the $q$-Gevrey sequences provide such examples.
\end{remark}

For the sake of completeness we comment on a minor technical issue and mention that for $p=1$ in the first estimate in \cite[$(3.3)$]{modgrowthstrangeII} we need normalization, i.e. $\mu_1=\frac{M_1}{M_0}\ge 1=:\mu_0$, but this technical assumption has not been made there. However, in the further considerations in \cite{modgrowthstrangeII}, see $(3.6)$ and $(3.9)$, one does have this property and hence the full validity of $(3.3)$.

\subsection{Second statement}\label{secondmainresult}
When imposing some more (standard) growth and regularity properties on $\mathbf{M}$, then we extend the list of characterizations in Theorem \ref{firstcharprop} by involving the function $\omega_{\mathbf{M}}$.

\begin{theorem}\label{secondcharprop}
Let $\mathbf{M}=(M_j)_{j\in\NN}\in\RR_{>0}^{\NN}$ be a log-convex sequence such that $\lim_{j\rightarrow+\infty}(M_j)^{\frac{1}{j}}=+\infty$ and let $p\in[1,\infty)$ be fixed. Then the following are equivalent:
\begin{itemize}
\item[$(i)$] $\mathbf{M}$ satisfies \hyperlink{dc}{$(\on{dc})$} (i.e. \eqref{dc}):
$$\exists\;D\ge 1\;\forall\;j\in\mathbb{N}:\;\;\;M_{j+1}\le D^{j+1}M_j.$$

\item[$(ii)$] We get
\begin{equation}\label{newintcond}
\exists\;h'>h>0:\;\;\;\int_{0}^{+\infty}e^{p\omega_{\mathbf{M}}(ht)-p\omega_{\mathbf{M}}(h't)}dt<+\infty.
\end{equation}

\item[$(iii)$] We get
\begin{equation}\label{logabsorb}
\exists\;A\ge 1\;\exists\;B>0\;\forall\;t>0:\;\;\;\omega_{\mathbf{M}}(t)+\log(t)\le\omega_{\mathbf{M}}(At)+B.
\end{equation}
\end{itemize}
\end{theorem}
\emph{Note:} Again, the validity of this characterization is not depending on $p$ and, moreover, the parameters $h'>h>0$ in \eqref{newintcond} can be chosen to be the same as in \eqref{newseriescond}. And \eqref{newintcond} is tailor-made to treat weighted $L^p$-type norms for which the weights are expressed in terms of $\exp\circ\omega_{\mathbf{M}}$ and with parameters $h',h>0$. On the other hand, \eqref{logabsorb} is clearly sufficient to have \eqref{newintcond} and this observation is important in the mixed setting; see Section \ref{mixedsection}.

\demo{Proof}
$(i)\Leftrightarrow(ii)$ Indeed, we verify the equivalence between \eqref{newseriescond} and \eqref{newintcond} and the conclusion holds by $(i)\Leftrightarrow(iv)$ in Theorem \ref{firstcharprop}.

First, we follow the proof of \cite[Lemma 2.9]{weightedentireinclusion1} and involve the \emph{integral representation formula} \eqref{assointrepr} in order to estimate as follows for all $t>0$ and any $h'>h>0$ (fixed):
$$\omega_{\mathbf{M}}(h't)-\omega_{\mathbf{M}}(ht)=\int_{ht}^{h't}\frac{\Sigma_{\mathbf{M}}(u)}{u}du\ge\Sigma_{\mathbf{M}}(ht)\int_{ht}^{h't}\frac{1}{u}du=\Sigma_{\mathbf{M}}(ht)\log\left(\frac{h'}{h}\right),$$
and similarly
$$\omega_{\mathbf{M}}(h't)-\omega_{\mathbf{M}}(ht)\le\Sigma_{\mathbf{M}}(h't)\log\left(\frac{h'}{h}\right).$$
For $t=0$ one has the trivial equalities $0=0$. Summarizing, so far we have shown:
$$\forall\;h'>h>0:\;\;\;\int_0^{+\infty}\left(\frac{h}{h'}\right)^{p\Sigma_{\mathbf{M}}(h't)}dt\le\int_0^{+\infty}e^{p\omega_{\mathbf{M}}(ht)-p\omega_{\mathbf{M}}(h't)}dt\le\int_0^{+\infty}\left(\frac{h}{h'}\right)^{p\Sigma_{\mathbf{M}}(ht)}dt.$$
Next, we treat the third integral and the first one is then completely analogous:
\begin{align*}
&\int_0^{+\infty}\left(\frac{h}{h'}\right)^{p\Sigma_{\mathbf{M}}(ht)}dt=\frac{\mu_1}{h}+\sum_{j=1}^{+\infty}\int_{\mu_j/h}^{\mu_{j+1}/h}\left(\frac{h}{h'}\right)^{pj}dt=\frac{\mu_1}{h}+\sum_{j=1}^{+\infty}\left(\frac{h}{h'}\right)^{pj}(\mu_{j+1}/h-\mu_j/h),
\end{align*}
and this last expression is finite if and only if $\sum_{j=1}^{+\infty}\left(\frac{h}{h'}\right)^{pj}(\mu_{j+1}-\mu_j)<+\infty$. Note that in the above identity the indices $j$ such that $\mu_j=\mu_{j+1}$ yield the value $0$ both for the corresponding integral and for the summand in the series (cf. the definition of the counting function in \eqref{counting}). Altogether, this gives the desired equivalence with the same $h'>h>0$.\vspace{6pt}

$(i)\Leftrightarrow(iii)$ This is shown, even in the mixed setting, in \cite[Lemma 2]{nuclearglobal2} and for $(i)\Rightarrow(iii)$ see also \cite[Prop. 3.4]{Komatsu73}. For the sake of completeness we repeat the arguments:

$(i)\Rightarrow(iii)$ By assumption we get $\frac{M_0t^j}{M_j}\le D\frac{M_0(Dt)^j}{M_{j+1}}=\frac{1}{t}\frac{M_0(Dt)^{j+1}}{M_{j+1}}$ for all $t>0$ and $j\in\NN$ (the first estimate also holds for $t=0$). Applying $\log$ to this estimate and recalling the definition of the associated function \eqref{logabsorb} follows (with $A:=D$ and $B>0$ arbitrary).

$(iii)\Rightarrow(i)$ Conversely, we involve \eqref{Prop32Komatsu} and get for all $j\in\NN$
\begin{align*}
M_j&=M_0\sup_{t>0}\frac{t^j}{\exp(\omega_{\mathbf{M}}(t))}\ge M_0e^{-B}\sup_{t>0}\frac{t^{j+1}}{\exp(\omega_{\mathbf{M}}(At))}= M_0e^{-B}\frac{1}{A^{j+1}}\sup_{s>0}\frac{s^{j+1}}{\exp(\omega_{\mathbf{M}}(s))}
\\&
=e^{-B}\frac{1}{A^{j+1}}M_{j+1},
\end{align*}
which yields the conclusion.
\qed\enddemo

\subsection{Third statement}\label{thirdmainresult}
The next aim is to characterize \eqref{dcMandelbrojt} in terms of $\omega_{\mathbf{M}}$, $\Sigma_{\mathbf{M}}$, and hence to extend the above results. We start with a technical observation.

\begin{proposition}\label{thirdcharprop}
Let $\mathbf{M}=(M_j)_{j\in\NN}\in\RR_{>0}^{\NN}$ be a log-convex sequence such that $\mathbf{M}_{\iota}=\lim_{j\rightarrow+\infty}(M_j)^{\frac{1}{j}}=+\infty$. Consider the following assertions:
\begin{itemize}
\item[$(i)$] $\omega_{\mathbf{M}}$ satisfies
\begin{equation}\label{thirdcharpropequ}
\liminf_{t\rightarrow+\infty}\frac{\omega_{\mathbf{M}}(t)}{(\log(t))^2}>0.
\end{equation}

\item[$(ii)$] $\omega_{\mathbf{M}}$ satisfies
\begin{equation}\label{thirdcharpropequconv}
\limsup_{t\rightarrow+\infty}\frac{\omega_{\mathbf{M}}(t)}{(\log(t))^2}<+\infty.
\end{equation}

\item[$(iii)$] $\Sigma_{\mathbf{M}}$ satisfies
$$\liminf_{t\rightarrow+\infty}\frac{\Sigma_{\mathbf{M}}(t)}{\log(t)}>0.$$

\item[$(iv)$] $\Sigma_{\mathbf{M}}$ satisfies
$$\limsup_{t\rightarrow+\infty}\frac{\Sigma_{\mathbf{M}}(t)}{\log(t)}<+\infty.$$

\item[$(v)$] $\mathbf{M}$ satisfies \hyperlink{dc}{$(\on{dc})$}.
\end{itemize}
Then the equivalences $(i)\Leftrightarrow(iii)\Leftrightarrow(v)$ and $(ii)\Leftrightarrow(iv)$ hold.
\end{proposition}

\demo{Proof}
Let $s,t\in\RR$ be such that $\max\{-\log(\mu_1),\log(\mu_1)\}<s<t$ and so, in any case, $0<s<t$. Then, since by definition $\Sigma_{\mathbf{M}}$ is non-decreasing, \eqref{assointrepr} implies
$$\omega_{\mathbf{M}}(e^t)=\int_{\mu_1}^{e^t}\frac{\Sigma_{\mathbf{M}}(\lambda)}{\lambda}d\lambda\le\Sigma_{\mathbf{M}}(e^t)\int_{\mu_1}^{e^t}\frac{1}{\lambda}d\lambda=\Sigma_{\mathbf{M}}(e^t)(t-\log(\mu_1))\le\Sigma_{\mathbf{M}}(e^t)2t,$$
and analogously (since $e^s>\mu_1$ holds):
$$\omega_{\mathbf{M}}(e^{t})=\int_{\mu_1}^{e^{t}}\frac{\Sigma_{\mathbf{M}}(\lambda)}{\lambda}d\lambda\ge\int_{e^{s}}^{e^{t}}\frac{\Sigma_{\mathbf{M}}(\lambda)}{\lambda}d\lambda\ge\Sigma_{\mathbf{M}}(e^{s})(t-s).$$
We choose $t:=2s(>s)$ and this gives
$$\forall\;s>\max\{-\log(\mu_1),\log(\mu_1)\}:\;\;\;\frac{\Sigma_{\mathbf{M}}(e^{s})}{s}\le\frac{\omega_{\mathbf{M}}(e^{2s})}{s^2}\le 8\frac{\Sigma_{\mathbf{M}}(e^{2s})}{2s},$$
equivalently
$$\forall\;u>\max\{\mu_1^{-1},\mu_1\}(\ge1):\;\;\;\frac{\Sigma_{\mathbf{M}}(u)}{\log(u)}\le\frac{\omega_{\mathbf{M}}(u^2)}{(\log(u))^2}\le 8\frac{\Sigma_{\mathbf{M}}(u^2)}{\log(u^2)}.$$
When setting $u^2:=v$ in the second quotient we infer $\frac{\omega_{\mathbf{M}}(u^2)}{(\log(u))^2}=4\frac{\omega_{\mathbf{M}}(v)}{(\log(v))^2}$, and thus the equivalences $(i)\Leftrightarrow(iii)$ and $(ii)\Leftrightarrow(iv)$ are established.\vspace{6pt}

$(v)\Rightarrow(iii)$ First, by $(i)\Leftrightarrow(ii)$ in Theorem \ref{firstcharprop} we get \eqref{mudc} and note that in this condition one has $D>1$ since $D\in(0,1]$ contradicts the assumption $\lim_{j\rightarrow+\infty}\mu_j=+\infty$. Second, put $\mathbf{M}^{+1}:=(M_{j+1})_{j\in\NN}$ with the corresponding quotient sequence $\mu^{+1}:=(\mu_{j+1})_{j\in\NN}$ (forward shifted sequence). Indeed, by taking into account this shift, we get the following properties:

$\mathbf{M}^{+1}$ is log-convex resp. satisfies \hyperlink{dc}{$(\on{dc})$} if and only if $\mathbf{M}$ does so. Moreover, \hyperlink{dc}{$(\on{dc})$} (for $\mathbf{M}$) amounts to have $\mathbf{M}^{+1}\hyperlink{preceq}{\preceq}\mathbf{M}$. And when $\lim_{j\rightarrow+\infty}\mu_j=+\infty$ then, in particular, $\mu_j\ge 1$ for all $j$ sufficiently large and hence $\mathbf{M}$ is eventually increasing. In this case \hyperlink{dc}{$(\on{dc})$} does precisely mean that $\mathbf{M}^{+1}$ and $\mathbf{M}$ are equivalent. Finally, recall the definition \eqref{counting}
$$\Sigma_{\mathbf{M}^{+1}}(t)=|\{j\in\NN_{>0}: \mu_{j+1}=\mu^{+1}_j\le t\}|,\;\;\;t\ge 0.$$
Thus $\Sigma_{\mathbf{M}^{+1}}(t)=\Sigma_{\mathbf{M}}(t)-1$ for all $t\ge\mu^{+1}_1=\mu_2$ and $\Sigma_{\mathbf{M}^{+1}}(t)=0$ for all $t\in[0,\mu^{+1}_1)=[0,\mu_2)$, hence $\Sigma_{\mathbf{M}^{+1}}(t)=\Sigma_{\mathbf{M}}(t)-1$ for all $t\ge\mu_1$ and $\Sigma_{\mathbf{M}^{+1}}(t)=0=\Sigma_{\mathbf{M}}(t)$ for all $t\in[0,\mu_1)$.

Accordingly, for any $D>1$ (fixed) we put
\begin{equation}\label{countingD}
\Sigma_D(t):=|\{j\in\NN_{>0}: D^{j+1}\le t\}|,\;\;\;t\ge 0,
\end{equation}
and note that assumption \eqref{mudc} and the definition of the counting functions immediately give $\Sigma_{\mathbf{M}^{+1}}(t)\ge\Sigma_D(t)$ for all $t\ge 0$ and hence by the above $\Sigma_{\mathbf{M}}(t)-1\ge\Sigma_D(t)$ for all $t\ge\mu_1$. Also note that $\Sigma_D(t)=0$ for all $t\in[0,D^2)$ (again since the empty set has cardinality $0$). Moreover, if $t\ge 0$ is such that $D^{j+1}\le t<D^{j+2}$ for some $j\in\NN_{>0}$, then $\Sigma_D(t)=j$ and $j\log(D)\le(j+1)\log(D)\le\log(t)$, $\log(t)\le(j+2)\log(D)\le 2j\log(D)$ and for all $j\ge 2$. Hence, $\frac{\log(t)}{2\log(D)}\le\Sigma_D(t)\le\frac{\log(t)}{\log(D)}$ for all $t\ge D^3$ and, summarizing, it is shown that
$$\forall\;t\ge D^3(\ge\mu_3\ge\mu_1):\;\;\;\Sigma_{\mathbf{M}}(t)-1\ge\Sigma_D(t)\ge\frac{\log(t)}{2\log(D)},$$
which obviously implies $\liminf_{t\rightarrow+\infty}\frac{\Sigma_{\mathbf{M}}(t)}{\log(t)}>0$.\vspace{6pt}

$(iii)\Rightarrow(v)$ By assumption we get
$$\exists\;\epsilon>0\;\exists\;t_{\epsilon}>1\;\forall\;t\ge t_{\epsilon}:\;\;\;\Sigma_{\mathbf{M}}(t)>\epsilon\log(t)=\log(t^{\epsilon}).$$
Since $\lim_{j\rightarrow+\infty}\mu_j=+\infty$ there exists $j_{\epsilon}\in\NN_{>0}$ such that $\mu_j\ge t_{\epsilon}$ for all $j\ge j_{\epsilon}$ and for all $t\ge\mu_{j_{\epsilon}}$ we can find $j\ge j_{\epsilon}$ such that $\mu_j\le t<\mu_{j+1}$. W.l.o.g. we can assume that $\mu_{j_{\epsilon}}\le t<\mu_{j_{\epsilon}+1}$ and then we take $t=\mu_j$ and get
$$\log(\mu_j^{\epsilon})<\Sigma_{\mathbf{M}}(\mu_j)=j.$$
Thus $\mu_j<(\exp(\epsilon^{-1}))^j$ for all $j\ge j_{\epsilon}$ satisfying $\mu_j<\mu_{j+1}$. Moreover, when $\mu_{j-1}<\mu_j=\dots=\mu_{j+i}$ for some $i\in\NN_{>0}$ and $j>j_{\epsilon}$, then the above argument gives $t<(\exp(\epsilon^{-1}))^{j-1}$ for all $t\in[\mu_{j-1},\mu_j)$ and so $\mu_j\le(\exp(\epsilon^{-1}))^{j-1}$ which implies $\mu_{j+k}\le(\exp(\epsilon^{-1}))^{j+k}$ for all $0\le k\le i$, too.

Summarizing, this shows $\mu_j\le(\exp(\epsilon^{-1}))^j$ for all $j\ge j_{\epsilon}$ and thus, when replacing $\epsilon$ by some $\epsilon_1\in(0,\epsilon)$ if necessary, we get $\mu_j\le(\exp(\epsilon_1^{-1}))^j$ for all $j\in\NN_{>0}$; i.e. \eqref{mudc}.
\qed\enddemo

In order to proceed, we have to recall and introduce some notation: For any $\ell>1$ define the sequence
\begin{equation}\label{qGevrey}
\mathbf{M}^{(\ell)}:=(M^{(\ell)}_j)_{j\in\NN},\hspace{15pt}M^{(\ell)}_j:=\ell^{j^2}.
\end{equation}
In the literature, $\mathbf{M}^{(\ell)}$ is also called a \emph{$q$-Gevrey sequence} (with parameter $\ell>1$) and the corresponding sequence of quotients is given by $\mu^{(\ell)}_j=\ell^{2j-1}$ which should be compared with the proof of $(v)\Rightarrow(iii)$ in Proposition \ref{thirdcharprop}, in particular see \eqref{countingD}.

By using the order relation \hyperlink{preceq}{$\preceq$} we see that \eqref{dcMandelbrojt} for a given sequence $\mathbf{M}\in\RR_{>0}^{\NN}$ precisely means that
\begin{equation}\label{dcMandelbrojtqGevrey}
\exists\;\ell>1:\;\;\;\mathbf{M}\hyperlink{preceq}{\preceq}\mathbf{M}^{(\ell)}.
\end{equation}
Next, for functions $\sigma,\tau:[0,+\infty)\rightarrow[0,+\infty)$ such that $\lim_{t\rightarrow+\infty}\sigma(t)=+\infty=\lim_{t\rightarrow+\infty}\tau(t)$ we write $\sigma\hypertarget{ompreceq}{\preceq}\tau$ if
\begin{equation}\label{bigOrelation}
\tau(t)=O(\sigma(t))\;\text{as}\;t\rightarrow+\infty.	
\end{equation}
Moreover, write $\sigma\hypertarget{sim}{\sim}\tau$, if $\sigma\hyperlink{ompreceq}{\preceq}\tau$ and $\tau\hyperlink{ompreceq}{\preceq}\sigma$ and in this case $\sigma$ and $\tau$ are called \emph{equivalent} (cf. \cite[Sect. 2.3]{ultradifferentiablecomparison}). Relation \hyperlink{ompreceq}{$\preceq$} is crucial in the ultradifferentiable weight function setting; we refer to \cite[Sect. 5.1, Lemma 5.16 \& Cor. 5.17]{compositionpaper}. Define $\log^2: t\mapsto(\log(t))^2$ and, consequently, \eqref{thirdcharpropequ} precisely means that $\log^2(t)=O(\omega_{\mathbf{M}}(t))$ as $t\rightarrow+\infty$; i.e. $\omega_{\mathbf{M}}\hyperlink{ompreceq}{\preceq}\log^2$. (Analogously \eqref{thirdcharpropequconv} means $\log^2\hyperlink{ompreceq}{\preceq}\omega_{\mathbf{M}}$, but this condition will not be relevant in our arguments.) Moreover, with $\omega_2(t):=\max\{0,(\log(t))^2\}$ it holds that $\omega_2\hyperlink{sim}{\sim}\log^2\hyperlink{sim}{\sim}t\mapsto(\log(1+t))^2$ and $\omega_2$ is even \emph{normalized;} i.e. $\omega_2(t)=0$ for all $t\in[0,1]$. We recall and summarize some known facts:

\begin{itemize}
\item[(1)] These functions, and more generally when $2$ is replaced by an arbitrary parameter $s>1$, play a prominent role as standard examples in the ultradifferentiable weight function setting in the sense of Braun-Meise-Taylor (cf. \cite{BraunMeiseTaylor90}) and we refer to \cite{BonetBraunMeiseTaylorWhitneyextension}, \cite{BonetMeiseMelikhov07}, \cite{dissertation}, \cite{compositionpaper}. Note that in \cite[Ex. 20]{BonetMeiseMelikhov07} it has been commented that these weights serve as weight functions such that the corresponding ultradifferentiable classes cannot be described alternatively by a weight sequence.

\item[(2)] Based on and inspired by the results from \cite{BonetMeiseMelikhov07}, in \cite{dissertation} and \cite{compositionpaper} a central idea has been to associate with each (normalized) BMT-weight function $\omega$ a \emph{weight matrix} $\mathcal{M}_{\omega}:=\{\mathbf{W}^{(\ell)}=(W^{(\ell)}_j)_{j\in\NN}: \ell>0\}$ by
	$$W^{(\ell)}_j:=\exp\left(\frac{1}{\ell}\varphi^{*}_{\omega}(\ell j)\right),$$
with $\varphi^{*}_{\omega}(x):=\sup\{x y-\varphi_{\omega}(y): y\ge 0\}$, $x\ge 0$, and $\varphi_{\omega}:=\omega\circ\exp$. It turns out that basically for this definition it is sufficient to deal with non-decreasing $\omega$ satisfying $\log(t)=o(\omega(t))$ as $t\rightarrow+\infty$ \emph{(``matrix admissible weights'');} see \cite[Sect. 2.9]{ultradifferentiablecomparison} for more details and see \cite[Def. 1.4, Sect. 2]{ultradifferentiablecomparison} for the definition and detailed comments concerning BMT-weight functions.

\item[(3)] This method is convenient to describe both the weight sequence and weight function setting in a uniform way and to transfer proofs. For concrete examples the computation of the precise expression of $\mathbf{W}^{(\ell)}$ in terms of $\omega$ might be complicated and cumbersome. But for the weights $\omega_s:=\max\{0,(\log(t))^s\}$, $s>1$, it is possible to derive by elementary calculations a closed expression and this study has also been motivated by \cite[Ex. 20]{BonetMeiseMelikhov07}.

    From now on $\mathbf{W}^{(\ell)}$ shall denote the sequences associated with the weight function $\omega_2$.

    In \cite[Sect. 3.10, p. 29-32; $(3.10.1)$]{dissertation} and \cite[Lemma 7.4, $(7.2)$]{convenientweightmatrix} (see also \cite[Sect. 5.5]{whitneyextensionweightmatrix}) we have obtained the following formula for $\mathbf{W}^{(\ell)}$:
    \begin{equation}\label{qGevreytransform}
    \forall\;\ell>0\;\forall\;j\in\NN:\;\;\;W^{(\ell)}_j=(\exp(\ell/4))^{j^2}=M^{(\exp(\ell/4))}_j;
    \end{equation}
    with other words $\mathbf{W}^{(\ell)}=\mathbf{M}^{(\exp(\ell/4))}$ for all $\ell>0$ and thus the matrix $\mathcal{M}_{\omega_2}$ precisely consists of all $q$-Gevrey sequences; i.e. as sets $\mathcal{M}_{\omega_2}=\mathcal{G}^q:=\{\mathbf{M}^{(\ell)}: \ell>0\}$. (One has $R(2)=\frac{1}{4}$ with $R(s):=\frac{1}{s^{\frac{1}{s-1}}}-\frac{1}{s^{\frac{s}{s-1}}}$, so $R(s)\rightarrow 1$ as $s\rightarrow+\infty$ but note that the second term in line 4 before \cite[$(3.10.1)$, p. 30]{dissertation} is incorrect and $\frac{s^{\frac{1}{s}}}{s}$ should be replaced by $s^{\frac{s}{1-s}}$ and which tends to $0$ as $s\rightarrow+\infty$.)

\item[(4)] Finally, we mention that for any given $\mathbf{M}\in\hyperlink{LCset}{\mathcal{LC}}$ we can associate the matrix $\mathcal{M}_{\omega_{\mathbf{M}}}$, see e.g. \cite[Lemma 2.8]{testfunctioncharacterization} resp. \cite[Lemma 2.4]{sectorialextensions} and the references mentioned in the proofs there. Moreover, in this case $\mathbf{M}=\mathbf{W}^{(1)}$, i.e. $\mathbf{M}$ coincides with the associated weight sequence with index $\ell=1$; see e.g. \cite[$(2.13)$, p. 407]{subaddlike}.
\end{itemize}

We are now ready to prove the next characterization for which we involve weight matrix techniques.

\begin{theorem}\label{fourthcharprop}
Let $\mathbf{M}=(M_j)_{j\in\NN}\in\RR_{>0}^{\NN}$ be a log-convex sequence such that $\mathbf{M}_{\iota}=\lim_{j\rightarrow+\infty}(M_j)^{\frac{1}{j}}=+\infty$. Then the following are equivalent:
\begin{itemize}
\item[$(i)$] $\mathbf{M}$ satisfies \eqref{dcMandelbrojt}.

\item[$(ii)$] $\omega_{\mathbf{M}}$ satisfies \eqref{thirdcharpropequ}.
\end{itemize}
\end{theorem}

\demo{Proof}
First, there exists $\mathbf{N}\in\hyperlink{LCset}{\mathcal{LC}}$ such that $M_j=N_j$ for all sufficiently large $j$ and hence
\begin{equation}\label{MNnormequiv}
\exists\;A\ge 1\;\forall\;j\in\NN:\;\;\;\frac{1}{A}N_j\le M_j\le AM_j.
\end{equation}
More details concerning the precise definition of $\mathbf{N}$ can be found in \cite[Rem. 2.5]{weightedentireinclusion1} and \cite[Prop. 3.5]{regularnew}.\vspace{6pt}

$(i)\Rightarrow(ii)$ \eqref{MNnormequiv} implies that $\mathbf{N}$ satisfies \eqref{dcMandelbrojt} too since this condition is obviously preserved under equivalence of sequences. Next, as mentioned above \eqref{dcMandelbrojt} precisely means (cf. \eqref{dcMandelbrojtqGevrey}):
$$\exists\;\ell_0>1:\;\;\;\mathbf{N}\hyperlink{preceq}{\preceq}\mathbf{M}^{(\ell_0)}.$$
Thus, since $N_0=1=M^{(\ell)}_0$, this yields
$$\exists\;\ell_0>1\;\exists\;C\ge 1\;\forall\;j\in\NN:\;\;\;N_j\le C^jM^{(\ell_0)}_j(=C^j\ell_0^{j^2}),$$
and so $\frac{t^j}{M^{(\ell_0)}_j}\le\frac{(Ct)^j}{N_j}$ for all $t\ge 0$ and $j\in\NN$. By definition of associated weight functions and \eqref{qGevreytransform} this estimate precisely gives
\begin{equation}\label{MNnormequiv0}
\exists\;\ell_0>1\;\exists\;C\ge 1\;\forall\;t\ge 0:\;\;\;\omega_{\mathbf{W}^{(4\log(\ell_0))}}(t)=\omega_{\mathbf{M}^{(\ell_0)}}(t)\le\omega_{\mathbf{N}}(Ct).
\end{equation}
Similarly, \eqref{MNnormequiv} implies
\begin{equation}\label{MNnormequiv1}
\exists\;A\ge 1\;\forall\;t\ge 0:\;\;\;\omega_{\mathbf{N}}(t)\le\omega_{\mathbf{M}}(t)+\log(A).
\end{equation}
Next recall that we have $\omega_2\hyperlink{sim}{\sim}\omega_{\mathbf{W}^{(\ell)}}$ for each $\ell>0$, more precisely
\begin{equation}\label{goodequivalenceclassic}
\forall\;\ell>0\;\exists\,D_{\ell}>0\;\forall\;t\ge 0:\;\;\;\ell\omega_{\mathbf{W}^{(\ell)}}(t)\le\omega_2(t)\le 2\ell\omega_{\mathbf{W}^{(\ell)}}(t)+D_{\ell};
\end{equation}
for this apply \cite[Theorem 4.0.3, Lemma 5.1.3]{dissertation} resp. \cite[Lemma 5.7]{compositionpaper} and also \cite[Lemma 2.5]{sectorialextensions} to the function $\omega_2$. \eqref{goodequivalenceclassic} for $\ell:=4\log(\ell_0)$ and the above estimates \eqref{MNnormequiv0} and \eqref{MNnormequiv1} give
\begin{align*}
&\exists\;\ell_0>1\;\exists\;A,C,D\ge 1\;\forall\;t\ge 0:
\\&
\frac{1}{8\log(\ell_0)}\omega_2(t/C)-\frac{D}{8\log(\ell_0)}\le\omega_{\mathbf{W}^{(4\log(\ell_0))}}(t/C)\le\omega_{\mathbf{N}}(t)\le\omega_{\mathbf{M}}(t)+\log(A),
\end{align*}
with $D:=D_{4\log(\ell_0)}$. Obviously $\omega_2(2t)=O(\omega_2(t))$ as $t\rightarrow+\infty$ holds, and this is a standard condition in the theory of weight functions in the sense of Braun-Meise-Taylor (see e.g. \cite{BraunMeiseTaylor90}, \cite{BonetMeiseMelikhov07}, \cite{compositionpaper}). By iterating this condition we get $\omega_2(t)\le L\omega_2(t/C)+L$ for some $L\ge 1$ and all $t\ge 0$; see also \cite[Lemma 2.1 (I)]{ultradifferentiablecomparison}. Thus the previous estimate implies
$$\exists\;B\ge 1\;\forall\;t\ge 0:\;\;\;\omega_2(t)\le B\omega_{\mathbf{M}}(t)+B;$$
i.e. $\omega_{\mathbf{M}}\hyperlink{ompreceq}{\preceq}\omega_2$. Finally, since $\omega_2\hyperlink{sim}{\sim}\log^2$ this proves $\omega_{\mathbf{M}}\hyperlink{ompreceq}{\preceq}\log^2$ and which is precisely \eqref{thirdcharpropequ}.\vspace{6pt}

$(ii)\Rightarrow(i)$ By \eqref{thirdcharpropequ} and $\omega_2\hyperlink{sim}{\sim}\log^2$ we get $\omega_{\mathbf{M}}\hyperlink{ompreceq}{\preceq}\omega_2$ and by \eqref{MNnormequiv1} this relation implies $\omega_{\mathbf{N}}\hyperlink{ompreceq}{\preceq}\omega_2$. Then invoke (the proof of) \cite[Lemma 5.16 (1)]{compositionpaper} resp. see also \cite[Sect. 6, (III), $(6.6)$]{modgrowthstrangeII} and obtain the following estimate (the equalities follows by \eqref{qGevreytransform}):
$$\exists\;B\ge 1\;\forall\;\ell>0\;\exists\;C\ge 1\;\forall\;j\in\NN:\;\;\;N^{(\ell)}_j\le CW^{(B\ell)}_j=CM^{(\exp(B\ell/4))}_j=C(\exp(B\ell/4))^{j^2}.$$
Here, $\mathbf{N}^{(\ell)}$ shall denote the sequences of the associated weight matrix $\mathcal{M}_{\omega_{\mathbf{N}}}$ and since $\mathbf{N}\in\hyperlink{LCset}{\mathcal{LC}}$ we have that $\mathbf{N}=\mathbf{N}^{(1)}$; recall (4) above. So we choose $\ell:=1$ and this estimate yields \eqref{dcMandelbrojt} for $\mathbf{N}$. In view of \eqref{MNnormequiv} then also $\mathbf{M}$ has \eqref{dcMandelbrojt} and we are done.
\qed\enddemo

\subsection{Fourth statement}\label{fourthmainresult}
According to \cite[p. 669]{matsumoto} and \cite[p. 104]{compositionpaper} very abstractly we say that a weighted class of smooth functions and being defined in terms of a sequence $\mathbf{M}$, write $\mathcal{A}_{\mathbf{M}}$ for this space, satisfies \emph{derivation closedness} if the following property holds: For each $f\in\mathcal{A}_{\mathbf{M}}$ also all (partial) derivatives of $f$ belong to $\mathcal{A}_{\mathbf{M}}$. For applications, especially in the context of PDE's, this property is natural and highly desirable. In \cite{matsumoto} globally defined spaces (of Roumieu-type) have been treated (using the symbol/functor $\mathcal{B}$) and \cite{compositionpaper} focuses on locally defined classes both of Roumieu- and Beurling-type with the symbol $\mathcal{E}$. We do not provide exact definitions of these spaces since this work focuses on growth properties for given sequences and results can be applied to different contexts. For precise definitions we refer to \cite{Komatsu73} and \cite{compositionpaper}.

We recall the following known result; see \cite[Prop. 2.12 (6), p. 104]{compositionpaper}.

\begin{theorem}\label{spacedc}
Let $\mathbf{M}=(M_j)_{j\in\NN}\in\RR_{>0}^{\NN}$ be a log-convex sequence such that $\mathbf{M}_{\iota}=\lim_{j\rightarrow+\infty}(M_j)^{\frac{1}{j}}=+\infty$. Then the following are equivalent:
\begin{itemize}
\item[$(i)$] The spaces $\mathcal{E}_{\{\mathbf{M}\}}$ (Roumieu-type) and/or $\mathcal{E}_{(\mathbf{M})}$ are/is derivation closed.

\item[$(ii)$] $\mathbf{M}$ satisfies \eqref{dc}.
\end{itemize}
\end{theorem}

\demo{Proof}
The argument uses the more general characterization of the inclusion relations, see \cite[Prop. 2.12 (1) \& Prop. 4.6]{compositionpaper}, for $\mathbf{M}$ and $\mathbf{M}^{+1}$ (forward-shifted sequence; see the proof of $(v)\Rightarrow(iii)$ in Proposition \ref{thirdcharprop}). Note that the same characterization holds for the globally defined classes in \cite{matsumoto} by the analogous proof and, indeed, these assertions are equivalent to $\mathbf{M}^{+1}\hyperlink{approx}{\approx}\mathbf{M}$.
\qed\enddemo

\subsection{Final result}
Gathering the information from the previous sections we arrive at the final and main characterizing result concerning derivation closedness; for the sake of convenience we list all assertions again and group them:

\begin{theorem}\label{maindcthm}
Let $\mathbf{M}=(M_j)_{j\in\NN}\in\RR_{>0}^{\NN}$ be a log-convex sequence such that $\mathbf{M}_{\iota}=\lim_{j\rightarrow+\infty}(M_j)^{\frac{1}{j}}=+\infty$ and let $p\in[1,+\infty)$. Then the following are equivalent:
\begin{itemize}
\item[$(i)$] The spaces $\mathcal{E}_{\{\mathbf{M}\}}$ (Roumieu-type) and/or $\mathcal{E}_{(\mathbf{M})}$ are/is derivation closed.

\item[$(ii)$] $\mathbf{M}$ satisfies
$$\exists\;D\ge 1\;\forall\;j\in\mathbb{N}:\;\;\;M_{j+1}\le D^{j+1}M_j.$$

\item[$(iii)$] $\mathbf{M}$ satisfies
$$\exists\;C\ge 1\;\forall\;j\in\NN:\;\;\;M_j\le C^{j^2+1}.$$

\item[$(iv)$] $\mathbf{M}$ resp. the corresponding sequence of roots satisfies
$$\exists\;A\ge 1\;\forall\;j\in\NN_{>0}:\;\;\;(M_{j+1})^{\frac{1}{j+1}}\le A(M_{j})^{\frac{1}{j}}.$$

\item[$(v)$] The sequence of quotients $\mu=(\mu_j)_{j\in\NN_{>0}}$ satisfies
$$\exists\;D\ge 1\;\forall\;j\in\mathbb{N}:\;\;\;\mu_{j+1}\le D^{j+1}.$$

\item[$(vi)$] The sequence of quotients $\mu=(\mu_j)_{j\in\NN_{>0}}$ satisfies
$$\exists\;A\ge 1\;\forall\;j\in\mathbb{N}_{>0}:\;\;\;\mu_{j+1}-\mu_j\le A^j.$$

\item[$(vii)$] The sequence of quotients $\mu=(\mu_j)_{j\in\NN_{>0}}$ satisfies
$$\exists\;h'>h>0:\;\;\;\sum_{j=1}^{+\infty}\left(\frac{h}{h'}\right)^{pj}(\mu_{j+1}-\mu_j)<+\infty.$$

\item[$(viii)$] The counting function $\Sigma_{\mathbf{M}}$ satisfies
$$\liminf_{t\rightarrow+\infty}\frac{\Sigma_{\mathbf{M}}(t)}{\log(t)}>0.$$

\item[$(ix)$] The associated function $\omega_{\mathbf{M}}$ satisfies
$$\liminf_{t\rightarrow+\infty}\frac{\omega_{\mathbf{M}}(t)}{(\log(t))^2}>0.$$

\item[$(x)$] The associated function $\omega_{\mathbf{M}}$ satisfies
$$\exists\;h'>h>0:\;\;\;\int_{0}^{+\infty}e^{p\omega_{\mathbf{M}}(ht)-p\omega_{\mathbf{M}}(h't)}dt<+\infty.$$

\item[$(xi)$] The associated function $\omega_{\mathbf{M}}$ satisfies
$$\exists\;A\ge 1\;\exists\;B>0\;\forall\;t>0:\;\;\;\omega_{\mathbf{M}}(t)+\log(t)\le\omega_{\mathbf{M}}(At)+B.$$
\end{itemize}
All listed assertions are preserved under equivalence of log-convex sequences.
\end{theorem}

\section{A (counter-)example}\label{counterexsection}
We construct a (counter-)example $\mathbf{M}$ which shows that without having log-convexity one cannot expect the full validity of Theorem \ref{maindcthm}. More precisely, we show that \eqref{dcMandelbrojt}, i.e. assertion $(vi)$ in Theorem \ref{firstcharprop} resp. $(iii)$ in Theorem \ref{maindcthm} is only necessary for having derivation closedness. On the other hand for the constructed sequence we have that $j\mapsto(M_j)^{\frac{1}{j}}$ is non-decreasing and hence still some regularity. Before stating and proving the main result we give the following definition:

\begin{definition}\label{nqadef}
A sequence $\mathbf{M}=(M_j)_{j\in\NN}\in\RR_{>0}^{\NN}$ such that $j\mapsto(M_j)^{\frac{1}{j}}$ is non-decreasing is called \emph{non-quasianalytic} if
\begin{equation}\label{nqadefequ}
\sum_{j\in\NN_{>0}}\frac{1}{(M_j)^{\frac{1}{j}}}<+\infty.
\end{equation}
If this series is divergent, then $\mathbf{M}$ is called \emph{quasianalytic.}
\end{definition}

This terminology is suggested by the general version of the \emph{Denjoy-Carleman theorem} by H\"{o}rmander, see \cite[Thm. 1.3.8]{hoermander}, for which log-convexity is not required necessarily and in the crucial series one deals with two regularizations of $\mathbf{M}$ instead: Consider, on the one hand, the quotients of the aforementioned log-convex minorant $\mathbf{M}^{\on{lc}}$ and, on the other hand, the increasing minorant of the sequence of roots $(M_j^{\frac{1}{j}})_{j\in\NN_{>0}}$. Thus the corresponding ultradifferentiable class does contain non-trivial functions with compact support if and only if the series of the increasing minorant of the sequence of roots is convergent; i.e. replace in \eqref{nqadefequ} the sequence of roots by its increasing minorant. Hence if $j\mapsto(M_j)^{\frac{1}{j}}$ is non-decreasing, as in the assumption in Definition \ref{nqadef}, then the corresponding ultradifferentiable class contains non-trivial functions with compact support if and only if \eqref{nqadefequ}.

\begin{proposition}\label{Mcounterexample}
There exist sequences $\mathbf{M}\in\RR_{>0}^{\NN}$ satisfying the following properties:
\begin{itemize}
\item[$(a)$] $\mathbf{M}$ is normalized (i.e. $1=M_0\le M_1$), $j\mapsto(M_j)^{\frac{1}{j}}$ is non-decreasing and $\mathbf{M}_{\iota}=\lim_{j\rightarrow+\infty}(M_j)^{\frac{1}{j}}=+\infty$, but $\mathbf{M}$ is not log-convex (and hence, in particular, $\mathbf{M}\notin\hyperlink{LCset}{\mathcal{LC}}$).

\item[$(b)$] $\mathbf{M}$ does not satisfy derivation closedness, i.e. \eqref{dc} fails, hence assertions $(i)-(iii)$ and $(v)$ in Theorem \ref{firstcharprop} fail too and so, in particular, \eqref{rootaincr} is also violated.

\item[$(c)$] \eqref{dcMandelbrojt} is valid, indeed even
\begin{equation}\label{dcMandelbrojtviol1}
\lim_{j\rightarrow+\infty}\frac{\log(M_j)}{j^2}=0\Leftrightarrow\lim_{j\rightarrow+\infty}\frac{\log((M_j)^{\frac{1}{j}})}{j}=0.
\end{equation}

\item[$(d)$] $(a)-(c)$ are not depending on the notion of \emph{(non-)quasianalyticity;} i.e. there exist both non-quasianalytic and quasianalytic $\mathbf{M}$ satisfying $(a)-(c)$.

\item[$(e)$] Any sequence $\mathbf{N}\in\RR_{>0}^{\NN}$ being equivalent to $\mathbf{M}$ satisfies $\lim_{j\rightarrow+\infty}(N_j)^{\frac{1}{j}}=+\infty$, $(b)$, $(c)$, \eqref{nqadefequ} if and only if $\mathbf{M}$ does so (cf. $(d)$). Finally, $\mathbf{N}$ cannot be log-convex and hence the whole equivalence class does not contain any log-convex weight sequence.
\end{itemize}
\end{proposition}

\emph{Note:} $(b)$ and $(c)$ together yield that in general the implications $(v)\Rightarrow(vi)$ resp. $(i)\Rightarrow(vi)$ in Theorem \ref{firstcharprop} cannot be reversed. The last statement in $(e)$ is in the spirit of \cite[Cor. 4.1]{regularnew} but the situation here is different and cannot be obtained by this result directly: necessarily we have in this present result $\lim_{j\rightarrow+\infty}(M_j)^{\frac{1}{j}}=+\infty$ and hence $\mathbf{M}_{\iota}=+\infty=\mathbf{M}_{\sigma}$ is valid.

\demo{Proof}
We construct $\mathbf{M}=(M_p)_{p\in\NN}$ in terms of auxiliary sequences. First, let $\mathbf{c}=(c_j)_{j\in\NN_{>0}}$, $\mathbf{d}=(d_j)_{j\in\NN_{>0}}$ be sequences of positive real numbers such that
\begin{itemize}
\item[$(i)$] $c_1=1=d_1$ (normalization);

\item[$(ii)$] $c_j<c_{j+1}$ and $d_j<d_{j+1}$ for all $j\in\NN_{>0}$ (both sequences are strictly increasing);

\item[$(iii)$] $\lim_{j\rightarrow+\infty}c_j=+\infty$ (unboundedness of $\mathbf{c}$);

\item[$(iv)$] $\sup_{j\in\NN_{>0}}d_{j+1}-d_j=+\infty$ (strong growth of $\mathbf{d}$).
\end{itemize}
By $(iv)$ we also get $\lim_{j\rightarrow+\infty}d_j=+\infty$. Then consider sequences of positive real numbers $\mathbf{a}=(a_j)_{j\in\NN_{>0}}$, $\mathbf{b}=(b_j)_{j\in\NN_{>0}}$ which are related to $\mathbf{c}$ and $\mathbf{d}$ as follows:
\begin{equation}\label{sequenceb}
b_j:=\exp(d_j),\;\;\;j\in\NN_{>0},
\end{equation}
\begin{equation}\label{sequencea}
a_1:=1,\hspace{15pt}a_j\in\NN_{>0},\;\;\;a_j+1<a_{j+1},\;\;\;a_j\ge d_jc_j,\;\;\;j\in\NN_{>0}.
\end{equation}
Therefore also $\mathbf{b}$ is strictly increasing and $\lim_{j\rightarrow+\infty}b_j=+\infty=\lim_{j\rightarrow+\infty}a_j$. Now set
\begin{equation}\label{Mdef}
M_0:=1,\hspace{20pt}M_p:=b_j^p,\;\;\;a_j\le p\le a_{j+1}-1,\;j\in\NN_{>0}.
\end{equation}
\emph{Claim $(a)$:} Normalization is clear because $1=M_0\le M_1=M_{a_1}=b_1=\exp(d_1)=\exp(1)$. Since $(M_p)^{\frac{1}{p}}=b_j$ for all $a_j\le p\le a_{j+1}-1$ and all $j\in\NN_{>0}$, clearly $p\mapsto(M_p)^{\frac{1}{p}}$ is non-decreasing and $\lim_{p\rightarrow+\infty}(M_p)^{\frac{1}{p}}=+\infty$ (see \eqref{Mdef}, \eqref{sequenceb} and $(iv)$).

However, $\mathbf{M}$ is \emph{not log-convex} since the sequence of quotients is given as follows: First $\mu_1=\frac{M_1}{M_0}=b_1$ and $\frac{M_p}{M_{p-1}}=\mu_p=b_j$ when $a_j<p\le a_{j+1}-1$, $j\in\NN_{>0}$. Second, when $p=a_{j+1}$ for some $j\in\NN_{>0}$, then $\mu_p=\frac{b_{j+1}^p}{b_j^{p-1}}>b_j=\mu_{p-1}$ since $j\mapsto b_j$ is strictly increasing (cf. \eqref{sequenceb} and $(ii)$). But for all $p=a_{j+1}+1$, $j\in\NN_{>0}$ arbitrary, we get $\mu_p=b_{j+1}<\frac{b_{j+1}^p}{b_j^{p-1}}=\mu_{p-1}$. Note that the second and the third case appear for these choices for any $j\in\NN_{>0}$ by assumption $a_j+1<a_{j+1}$. Finally, in view of \eqref{sequenceb} and $(iv)$ we even get
$$\sup_{j\in\NN_{>0}}\frac{\mu_{a_{j+1}}}{\mu_{a_{j+1}+1}}=\sup_{j\in\NN_{>0}}\left(\frac{b_{j+1}}{b_j}\right)^{a_{j+1}-1}=\sup_{j\in\NN_{>0}}\exp((a_{j+1}-1)(d_{j+1}-d_j))=+\infty.$$\vspace{6pt}

\emph{Claim $(b)$:} \eqref{dc} fails. First, by $(i)\Leftrightarrow(v)$ in Theorem \ref{firstcharprop} and by part $(a)$, which yields $\mathbf{M}_{\iota}>0$, we have to show that \eqref{rootaincr} does not hold. For this choose $p:=a_{j+1}-1$, $j\in\NN_{>0}$, and then similarly as before we use \eqref{sequenceb} and $(iv)$ and get $$\sup_{p\in\NN_{>0}}\frac{(M_{p+1})^{\frac{1}{p+1}}}{(M_{p})^{\frac{1}{p}}}\ge\sup_{j\in\NN_{>0}}\frac{(M_{a_{j+1}})^{\frac{1}{a_{j+1}}}}{(M_{a_{j+1}-1})^{\frac{1}{a_{j+1}-1}}}=\sup_{j\in\NN_{>0}}\frac{b_{j+1}}{b_j}=\sup_{j\in\NN_{>0}}\exp(d_{j+1}-d_j)=+\infty.$$
\vspace{6pt}

\emph{Claim $(c)$:} \eqref{dcMandelbrojtviol1} is valid. We investigate the quotient $\frac{\log((M_p)^{\frac{1}{p}})}{p}$ and note that
$$\forall\;j\in\NN_{>0}\;\forall\;p\in\NN_{>0},\;a_j\le p\le a_{j+1}-1:\;\;\;\frac{\log((M_p)^{\frac{1}{p}})}{p}=\frac{\log(b_j)}{p}=\frac{d_j}{p},$$
hence for all such indices $p$ (when $j$ is arbitrary but fixed) the maximum of this quotient is attained at $p=a_j$ and which gives the value $\frac{d_j}{a_j}\le\frac{1}{c_j}$ (see \eqref{sequencea}). By $(iii)$ this expression tends to $0$ as $j\rightarrow+\infty$.

Note that $\mathbf{c}$ can be chosen to grow arbitrarily fast, i.e. the rate of decay resp. rate of tending to $0$ in \eqref{dcMandelbrojtviol1} can be chosen as fast as desired, when the definition of $\mathbf{M}$ is changed accordingly (more precisely the sequence $\mathbf{a}$ is adjusted via the definition in \eqref{sequencea}). The same comment concerning the rate of increase applies to the sequence $\mathbf{d}$ and hence here one has to modify $\mathbf{b}$ \emph{and} $\mathbf{a}$ accordingly.\vspace{6pt}

\emph{Claim $(d)$:} First, since $p\mapsto(M_p)^{\frac{1}{p}}$ is non-decreasing (cf. $(a)$) in view of Definition \ref{nqadef} and by \eqref{Mdef} it holds that for the crucial series under consideration we get
\begin{equation}\label{nqacondition}
\sum_{p\ge 1}\frac{1}{(M_p)^{\frac{1}{p}}}=\sum_{j\in\NN_{>0}}\frac{a_{j+1}-a_j}{b_j}=\sum_{j\in\NN_{>0}}\frac{a_{j+1}-a_j}{\exp(d_j)}.
\end{equation}
Thus, in order to ensure non-quasianalyticity for $\mathbf{M}$ it suffices to have
$$\exists\;j_0\in\NN_{>0}\;\forall\;j\ge j_0:\;\;\;\frac{a_{j+1}-a_j}{\exp(d_j)}\le\frac{1}{2^j},$$
and when we assume that both $\mathbf{d}$ and $\mathbf{c}$ are sequences of integers subject to the assumptions $(i)-(iv)$, then the idea is to set $a_j:=d_jc_j$ in \eqref{sequencea} (with $1=a_1=d_1c_1=1$). Hence, in this situation it suffices to have
\begin{equation}\label{nqasuffcond}
\exists\;j_0\in\NN_{>0}\;\forall\;j\ge j_0:\;\;\;2^jd_{j+1}c_{j+1}\le\exp(d_j)+d_jc_j.
\end{equation}
Let us check that the particular choices $c_j:=j$ and $d_j:=j!$ satisfy \eqref{nqasuffcond}. Indeed, $(i)-(iv)$ are clear and concerning $(iv)$ note that $(j+1)!-j!=j!(j+1)-j!=j!j$. Then
$$2^j(j+1)(j+1)!-jj!\le 2^j2j(j+1)!-jj!=2^{j+1}j^2j!\le j!^3\le\exp(j!)$$
is valid for all $j\in\NN_{>0}$ sufficiently large. Finally, note that $a_{j+1}-a_j=(j+1)(j+1)!-jj!=j!((j+1)^2-j)=j!(j^2+j+1)>1$ for all $j\in\NN_{>0}$.\vspace{6pt}

On the other hand, in order to ensure quasianalyticity in view of \eqref{nqacondition} it suffices to have
\begin{equation}\label{qasuffcond}
\exists\;\epsilon>0\;\exists\;(j_n)_{n\in\NN_{>0}}\;\forall\;n\in\NN_{>0}:\;\;\;\frac{a_{j_n+1}-a_{j_n}}{\exp(d_{j_n})}\ge\epsilon.
\end{equation}
In order to guarantee \eqref{qasuffcond} (with $\epsilon=1$ and $j_n=n$ for all $n$) fix $\mathbf{c}$ satisfying $(i)-(iii)$ and such that $c_j\in\NN_{>0}$ for all $j$. Then introduce iteratively $\mathbf{d}$ with $d_1:=1$ and when $d_j$ is defined, then select $d_{j+1}\in\NN_{>0}$ such that $d_{j+1}\ge j+d_j$, which ensures $(ii)$ and $(iv)$, and such that $d_{j+1}c_{j+1}\ge\exp(d_j)+c_jd_j$. Finally, set (again) $a_j:=d_jc_j$, hence this last estimate also ensures $a_{j+1}-a_j>1$ for all $j\in\NN_{>0}$, $b_j:=\exp(d_j)$ and $\mathbf{M}$ is given according to \eqref{Mdef}.\vspace{6pt}

\emph{Claim $(e)$:} Relation \hyperlink{approx}{$\approx$} obviously preserves \eqref{dc}, \eqref{dcMandelbrojtviol1} and \eqref{nqadefequ}. Thus, if $\mathbf{N}\in\RR_{>0}^{\NN}$ satisfies $\mathbf{N}\hyperlink{approx}{\approx}\mathbf{M}$, then for this sequence \eqref{dc} fails, \eqref{dcMandelbrojtviol1} holds, and it satisfies \eqref{nqadefequ} if and only if $\mathbf{M}$ does so. Moreover, $\lim_{j\rightarrow+\infty}(N_j)^{\frac{1}{j}}=+\infty$ follows immediately. However, in general it is not clear if $j\mapsto(N_j)^{\frac{1}{j}}$ is non-decreasing and if $\mathbf{N}$ is normalized.

In any case, we show that $\mathbf{N}$ cannot be log-convex: If $\mathbf{N}$ is log-convex, since $\lim_{j\rightarrow+\infty}(N_j)^{\frac{1}{j}}=+\infty$ is valid, we can apply Theorem \ref{firstcharprop} to $\mathbf{N}$ but then the equivalence $(i)\Leftrightarrow(v)$ there contradicts assertions $(b)$ and $(c)$ for $\mathbf{N}$.
\qed\enddemo

For the constructed sequence $\mathbf{M}$ in Proposition \ref{Mcounterexample} assertions $(a)-(c)$ and $(e)$ are still true for any positive power sequence $\mathbf{M}^{\alpha}:=((M_j)^{\alpha})_{j\in\NN}$, $\alpha>0$. Note that the sequence of quotients of $\mathbf{M}^{\alpha}$ is given by the power $\alpha$ of the sequence $(\mu_j)_{j\in\NN_{>0}}$. In general, \eqref{nqadefequ} is not preserved when taking powers of the sequence but if \eqref{nqadefequ} holds then also for all $\mathbf{M}^{\alpha}$, $\alpha>1$, and if \eqref{nqadefequ} fails then also for all $\mathbf{M}^{\alpha}$, $0<\alpha<1$. Finally, $\mathbf{M}^{\alpha}$ and $\mathbf{M}^{\beta}$ are not equivalent for any $\alpha\neq\beta$ since $$\forall\;\alpha>\beta>0:\;\;\;\lim_{j\rightarrow+\infty}\left(\frac{M^{\alpha}_j}{M^{\beta}_j}\right)^{\frac{1}{j}}=\lim_{j\rightarrow+\infty}\left((M_j)^{\frac{1}{j}}\right)^{\alpha-\beta}=+\infty,$$
because $\lim_{j\rightarrow+\infty}(M_j)^{\frac{1}{j}}=+\infty$. Therefore, we obtain:

\begin{corollary}\label{Mcounterexamplecor}
There exist uncountable many pairwise non-equivalent sequences satisfying assertions $(a)-(e)$ in Proposition \ref{Mcounterexample}.
\end{corollary}

\section{Supplementary results and observations}\label{supplementarysection}

\subsection{On the stability of derivation closedness under modifications of sequences}\label{stabilitysection}
For applications it is natural to construct new sequences from given ones in a precise (algebraic) way and it is interesting to know if \hyperlink{dc}{$(\on{dc})$} is preserved under these actions; e.g. see the previous comments verifying Corollary \ref{Mcounterexamplecor} when involving powers of a given sequence and also the transformation in \cite[Sect. 3.1, (b), $(3.1)$]{modgrowthstrange}. So let $\mathbf{M}=(M_j)_{j\in\NN}\in\RR_{>0}^{\NN}$, $\mathbf{N}=(N_j)_{j\in\NN}\in\RR_{>0}^{\NN}$ be given.

\begin{itemize}
\item[$(a)$] Set $\mathbf{M}+\mathbf{N}:=(M_j+N_j)_{j\in\NN}$ \emph{(pointwise sum).} If both sequences $\mathbf{M}$ and $\mathbf{N}$ satisfy \hyperlink{dc}{$(\on{dc})$}, then also $\mathbf{M}+\mathbf{N}$. Since we have the trivial estimate $M_j+N_j\ge 2\min\{M_j,N_j\}$ for all $j$, the following holds: If either $\mathbf{M}_{\iota}>0$ or $\mathbf{N}_{\iota}>0$ then $(\mathbf{M}+\mathbf{N})_{\iota}>0$. If either $\mathbf{M}$ or $\mathbf{N}$ satisfies \eqref{rootinfinity}, then also  $\mathbf{M}+\mathbf{N}$.

\item[$(b)$] Set $\mathbf{M}\cdot\mathbf{N}:=(M_jN_j)_{j\in\NN}$ \emph{(pointwise product).} If both sequences $\mathbf{M}$ and $\mathbf{N}$ satisfy \hyperlink{dc}{$(\on{dc})$}, then also $\mathbf{M}\cdot\mathbf{N}$. If $\mathbf{M}_{\iota}>0$ and $\mathbf{N}_{\iota}>0$ resp. if both sequences satisfy \eqref{rootinfinity}, then also $\mathbf{M}\cdot\mathbf{N}$. If both sequences are log-convex, then $\mathbf{M}\cdot\mathbf{N}$ too and the sequence of quotients is obtained by taking the pointwise product of $\mu$ and $\nu$. Consequently, the basic assumptions and equivalent assertions listed in Theorem \ref{maindcthm} are preserved under the pointwise product.

\item[$(c)$] Set $\frac{\mathbf{M}}{\mathbf{N}}:=(\frac{M_j}{N_j})_{j\in\NN}$ \emph{(pointwise quotient).} If $\mathbf{M}$ satisfies \hyperlink{dc}{$(\on{dc})$} and $\mathbf{N}$ satisfies $N_j\le B^{j+1}N_{j+1}$ for some $B\ge 1$ and all $j\in\NN$, then also $\frac{\mathbf{M}}{\mathbf{N}}$ enjoys \hyperlink{dc}{$(\on{dc})$}. The condition on $\mathbf{N}$ is weak and clear if $\mathbf{N}$ is \emph{eventually increasing;} i.e. if there exists $j_0\in\NN_{>0}$ such that $N_j\le N_{j+1}$ for all $j\ge j_0$.

    Moreover, $\frac{\mathbf{M}}{\mathbf{N}}_{\iota}>0$ amounts to having $\mathbf{N}\hyperlink{preceq}{\preceq}\mathbf{M}$, whereas \eqref{rootinfinity} does precisely mean $\lim_{j\rightarrow+\infty}\left(\frac{M_j}{N_j}\right)^{\frac{1}{j}}=+\infty$. (This stronger relation is occasionally abbreviated in the literature by $\mathbf{N}\vartriangleleft\mathbf{M}$; see e.g. \cite[p. 101]{compositionpaper}.) However, log-convexity for $\frac{\mathbf{M}}{\mathbf{N}}$ is not clear in general even if both $\mathbf{M}$ and $\mathbf{N}$ are so; e.g. consider (different) Gevrey sequences $\mathbf{G}^s$, $s>0$. In order to transfer regularity properties to the quotient sequence it is unavoidable to assume that the growth of $\mathbf{M}$ and $\mathbf{N}$ is related in a precise way.

    Concerning $\mathbf{M}\cdot\mathbf{N}$ and $\frac{\mathbf{M}}{\mathbf{N}}$ we refer also to \cite{genLegendreconj}.

\item[$(d)$] For the next definition and properties we refer to \cite[Sect. 3.1]{convolvedmatrix} and the references there. Introduce the \emph{convolved sequence} $\mathbf{M}\star\mathbf{N}=(M\star N_j)_{j\in\NN}$ by
$$M\star N_{j}:=\min_{0\le k\le j}M_{k}N_{j-k},\;\,\;j\in\NN.$$
First, in view of \cite[Sect. 3.1 $(c)$]{convolvedmatrix}, i.e. $\mathbf{M}\star\mathbf{N}\le\max\{M_0,N_0\}\min\{\mathbf{M},\mathbf{N}\}$, it holds that when either $\mathbf{M}$ or $\mathbf{N}$ satisfies \hyperlink{dc}{$(\on{dc})$}, then also $\mathbf{M}\star\mathbf{N}$. It is not clear that in general $\mathbf{M}\star\mathbf{N}_{\iota}>0$ holds even if $\mathbf{M}_{\iota}>0$ and $\mathbf{N}_{\iota}>0$. On the other hand, by \cite[Sect. 3.1 $(e)$, Rem. 3.1]{convolvedmatrix} we have that if both $\mathbf{M}$ and $\mathbf{N}$ are log-convex, then $\min\{\mathbf{M}_{\iota},\mathbf{N}_{\iota}\}=\mathbf{M}\star\mathbf{N}_{\iota}$. In this situation also  $\mathbf{M}\star\mathbf{N}$ is log-convex which follows again by \cite[Sect. 3.1 $(e)$]{convolvedmatrix}. Therefore, the basic assumptions and equivalent assertions listed in Theorem \ref{maindcthm} are preserved under the convolution.
\end{itemize}

\subsection{On the mixed weight sequence setting}\label{mixedsection}
It is natural to ask if Theorem \ref{maindcthm} transfers to the mixed setting resp. which parts hold in this more general approach.
Most importantly, \eqref{dc} is replaced by
\begin{equation}\label{mixeddc}
\exists\;D\ge 1\;\forall\;j\in\NN:\;\;\;M_{j+1}\le D^{j+1}N_j,
\end{equation}
and let us abbreviate this relation by $(\mathbf{M},\mathbf{N})_{\on{dc}}$. Note: For any given $\mathbf{M}\in\RR_{>0}^{\NN}$ we can find $\mathbf{N}$ such that $(\mathbf{M},\mathbf{N})_{\on{dc}}$ by simply taking $\mathbf{N}\equiv\mathbf{M}^{+1}$ (forward-shifted sequence). \eqref{mixeddc} is natural when dealing with (abstractly given) weight matrices $\mathcal{N}:=\{\mathbf{N}^{(\ell)}: \ell>0\}$, see \cite[Sect. 4.1 \& 5]{compositionpaper} and this should be compared with the matrix $\mathcal{G}^q$ of $q$-Gevrey sequences mentioned above. Indeed, in this context due to the appearance of the matrix parameter one has $\mathbf{M}\equiv\mathbf{N}^{(\ell_1)}$, $\mathbf{N}\equiv\mathbf{N}^{(\ell_2)}$, and one deals with two kind of conditions: the ``Roumieu-type'' $(\mathcal{M}_{\{\on{dc}\}})$ ($\forall\;\ell_1>0\;\exists\;\ell_2>0$) and the ``Beurling-type''  $(\mathcal{M}_{(\on{dc})})$ ($\forall\;\ell_2>0\;\exists\;\ell_1>0$)) which allows to iterate \eqref{mixeddc}. Note that weight matrices $\mathcal{M}_{\omega}$ associated with (normalized) BMT-weight functions $\omega$ satisfy automatically both kind of conditions which follows by convexity of the conjugate $\varphi^{*}_{\omega}$ with $\mathbf{M}\equiv\mathbf{W}^{(\ell)}$, $\mathbf{N}\equiv\mathbf{W}^{(2\ell)}$ for arbitrary $\ell>0$; see \cite[$(5.6)$, Cor. 5.15]{compositionpaper}.

We summarize which parts of the results from Section \ref{mainresultsection} transfer to the mixed framework and comment on the failure of some implications as well:

\begin{itemize}
\item[$(A)$] Theorem \ref{spacedc} transfers immediately to the general weight matrix setting; see \cite[Prop. 4.6 (6)]{compositionpaper}.

\item[$(B)$] The equivalence between \eqref{mixeddc} and the mixed variant of \eqref{logabsorb} is already shown in \cite[Lemma 2]{nuclearglobal2}; normalization (i.e. $M_0=1=N_0$) is not effecting the proof. Again, iterating the mixed version of \eqref{logabsorb} implies the fact that one can absorb any multiple of the logarithm; see e.g. the proof of \cite[Lemma 4]{nuclearglobal2}. Since \eqref{logabsorb} is crucial when dealing with weighted $L^p$-type norms, $p\in[1,+\infty)$ also in the matrix setting one can get information.

\item[$(C)$] When iterating mixed conditions one has to be careful: If the number of iterations is fixed, like in \cite[Lemma 4]{nuclearglobal2} where it is only depending on the dimension or like in Remark \ref{mgcharrem} where the fixed parameter $d$ is relevant, then one can admit a controlled loss of regularity expressed in terms of (different) sequences resp. matrix parameters. However, the number of iterations in $(v)\Rightarrow(i)$, $(iii)\Rightarrow(ii)$ in the proof of Theorem \ref{firstcharprop} cannot be bounded uniformly and hence these steps are unclear.

\item[$(D)$] $(i)\Rightarrow(v)$ in Theorem \ref{firstcharprop} transfers immediately: If $(\mathbf{M},\mathbf{N})_{\on{dc}}$ and $\mathbf{N}_{\iota}>0$, then
$$\exists\;A\ge 1\;\forall\;j\in\NN_{>0}:\;\;\;(M_{j+1})^{\frac{1}{j+1}}\le A(N_{j})^{\frac{1}{j}};$$
i.e. the mixed variant of \eqref{rootaincr}.

\item[$(E)$] The proof of the equivalence $(iii)\Leftrightarrow(iv)$ in Theorem \ref{firstcharprop} transfers as follows: Assume that $\mathbf{M},\mathbf{N}\in\RR_{>0}^{\NN}$ such that $\mu_j\le\nu_{j+1}$ for all $j\in\NN_{>0}$, then
\begin{equation}\label{mudiffmix}
\exists\;A\ge 1\;\forall\;j\in\mathbb{N}_{>0}:\;\;\;\nu_{j+1}-\mu_j\le A^j
\end{equation}
if and only if
\begin{equation}\label{newseriescondmix}
\exists\;h'>h>0:\;\;\;\sum_{j=1}^{+\infty}\left(\frac{h}{h'}\right)^{pj}(\nu_{j+1}-\mu_j)<+\infty.
\end{equation}
Note that for the proof \eqref{newseriescondmix}$\Rightarrow$\eqref{mudiffmix} following $(iv)\Rightarrow(iii)$ in Theorem \ref{firstcharprop} the first estimate holds since by $\mu_j\le\nu_{j+1}$ each summand is non-negative.

\item[$(F)$] In general $(\mathbf{M},\mathbf{N})_{\on{dc}}$ does not imply \eqref{mudiffmix}; we verify this for the particular choice $\mathbf{N}\equiv\mathbf{M}^{+1}$. First, as mentioned before, obviously $(\mathbf{M},\mathbf{M}^{+1})_{\on{dc}}$ is valid and for $\mathbf{N}\equiv\mathbf{M}^{+1}$ we get $\nu_{j+1}=\mu^{+1}_{j+1}=\mu_{j+2}$. Define $\mathbf{M}$ in the spirit of \eqref{Mdef} as follows:
\begin{equation*}
M_0:=1,\hspace{20pt}\mu_p:=b_j,\;\;\;a_j\le p\le a_{j+1}-1,\;\;\;j\in\NN_{>0}.
\end{equation*}
Here $\mathbf{a}:=(a_j)_{j\in\NN_{>0}}$ is an arbitrary sequence of integers such that
$$a_1=1,\hspace{15pt}a_j+1<a_{j+1},\;j\in\NN_{>0},$$
and for $\mathbf{b}:=(b_j)_{j\in\NN_{>0}}$ assume
$$b_1=1,\hspace{15pt}b_j<b_{j+1},\hspace{15pt}b_{j+1}-b_j\ge j^j,\;\;\;j\in\NN_{>0}.$$
Then for $p\in\NN_{>0}$ put $M_p:=\prod_{k=1}^p\mu_k$ and by definition clearly $\mathbf{M}\in\hyperlink{LCset}{\mathcal{LC}}$. However, since $a_{j+1}-2\ge a_j$ for all $j\in\NN_{>0}$ it follows that
$$\forall\;j\in\NN_{>0}:\;\;\;\mu_{a_{j+1}}-\mu_{a_{j+1}-2}=b_{j+1}-b_j\ge j^j,$$
which obviously contradicts \eqref{mudiffmix} for $\nu\equiv\mu^{+1}$.

\item[$(G)$] The equivalence between \eqref{newseriescondmix} and the corresponding mixed variant of \eqref{newintcond}, i.e. $$\int_{0}^{+\infty}e^{p\omega_{\mathbf{N}}(ht)-p\omega_{\mathbf{M}}(h't)}dt<+\infty,$$
    is unclear: Following the given arguments in the above proof, for log-convex $\mathbf{M}$, $\mathbf{N}$ satisfying $\mathbf{M}_{\iota}=+\infty=\mathbf{N}_{\iota}$ and $\mu_j\le\nu_j$ for all $j\in\NN_{>0}$, for $\omega_{\mathbf{N}}(ht)-\omega_{\mathbf{M}}(h't)$ we can only prove an upper bound in terms of $\Sigma_{\mathbf{M}}(h't)\log(h'/h)$ resp. a lower bound in terms of $\Sigma_{\mathbf{N}}(ht)\log(h'/h)$. But these expressions refer to \eqref{newseriescond} for $\mathbf{M}$ resp. for $\mathbf{N}$ and hence to \eqref{dc} for $\mathbf{M}$ resp. for $\mathbf{N}$ and which fails in this setting in general.

    However, in view of comments $(B)$, $(E)$ and $(F)$ this observation is not too surprising: $(\mathbf{M},\mathbf{N})_{\on{dc}}$ guarantees always control on the mixed integral but \eqref{newseriescondmix} is violated in general.

\item[$(H)$] The generalizations of $(i)\Leftrightarrow(iii)\Leftrightarrow(v)$ in Proposition \ref{thirdcharprop} and Theorem \ref{fourthcharprop} are not relevant in this context in the following sense: One can have $(\mathbf{M},\mathbf{N})_{\on{dc}}$ but neither \eqref{dcMandelbrojt} for $\mathbf{M}$ and/or $\mathbf{N}$ nor \eqref{thirdcharpropequ} for $\omega_{\mathbf{M}}$ and/or $\omega_{\mathbf{N}}$ holds. We fix $\ell>0$ and take $\mathbf{M}\equiv\mathbf{W}^{(\ell)}$, $\mathbf{N}\equiv\mathbf{W}^{(2\ell)}$, for a (normalized) BMT-weight function $\omega$ satisfying $\omega(t)=o((\log(t))^2)$ as $t\rightarrow+\infty$; take e.g. $\omega_s:=\max\{0,(\log(t))^s\}$ with $s\in(1,2)$. By the above $(\mathbf{W}^{(\ell)},\mathbf{W}^{(2\ell)})_{\on{dc}}$ is true but the associated sequences violate \eqref{dcMandelbrojt}; see again the precise formulas in \cite[Sect. 3.10, p. 29-32; $(3.10.1)$]{dissertation}, \cite[Lemma 7.4, $(7.2)$]{convenientweightmatrix} and also \cite[Sect. 5.5]{whitneyextensionweightmatrix}. And in view of \eqref{goodequivalenceclassic} the corresponding associated functions violate \eqref{thirdcharpropequ}. (On the other hand, $(ii)\Leftrightarrow(iv)$ in Proposition \ref{thirdcharprop} is valid for these associated weight sequences.)
\end{itemize}

\bibliographystyle{plain}
\bibliography{Bibliography}
\end{document}